\documentclass{amsart}

\numberwithin{equation}{section}

\usepackage{amssymb}
\usepackage{enumerate}
\usepackage{xspace}
\usepackage{a4wide}
\usepackage{pstcol}

\newtheorem{thm}{Theorem}[section]
\newtheorem{lem}[thm]{Lemma}
\newtheorem{cor}[thm]{Corollary} 
\newtheorem{sublem}[thm]{Sub-lemma}
\newtheorem{prop}[thm]{Proposition}

\newtheorem{defin}{Definition}
\newtheorem{cond}{Condition}

\newtheorem{rem}[thm]{Remark}

\newcommand\B{{\mathcal B}}
\newcommand\Co{{\mathcal C}}
\newcommand\D{{\mathcal D}}

\newcommand\Lp{{\mathcal L}}
\newcommand\M{{\mathcal M}}
\newcommand\Or{{\it O}}
\newcommand\T{{\mathcal T}}
\newcommand\W{{\mathcal W}}

\newcommand\A{{\mathbb A}}
\newcommand\C{{\mathbb C}}

\newcommand\N{{\mathbb N}}
\newcommand\R{{\mathbb R}}

\newcommand\ve{\varepsilon}
\newcommand\vf{\varphi}

\newcommand\Id{\text{\bf Id}}

\newcommand{\co}{C}

\newcommand{\tp}{{\bf p}}
\newcommand{\bphi}{\bar\phi}

\begin{document}

\title{On Contact Anosov Flows}

\author{Carlangelo Liverani}
\address{Carlangelo Liverani\\
Dipartimento di Matematica\\
II Universit\`{a} di Roma (Tor Vergata)\\
Via della Ricerca Scientifica, 00133 Roma, Italy.}
\email{{\tt liverani@mat.uniroma2.it}}
\thanks{It is a pleasure to thank Lai-Sang Young for many discussions
on the subject without which this paper would not exist. I also
profited from several conversations with V.Baladi, D.Dolgopyat, F.Ledrappier and
S.Luzzatto. In addition, I thank M.Pollicott and the anonymous
referees for pointing out several
imprecisions in previous versions. I acknowledge
the partial support of the ESF Programme PRODYN and the hospitality of
Courant Istitute and I.H.E.S. where part of the paper was written.} 
\date{March 14, 2003}
\begin{abstract}
Exponential decay of correlations for $\Co^{4}$ Contact Anosov flows is
established. This implies, in particular, exponential decay of
correlations for all smooth geodesic flows in strictly negative curvature. 
\end{abstract}
\maketitle

\section{Introduction}

The study of decay of correlations for hyperbolic systems
goes back to the work of Sinai \cite{Si0} and Ruelle \cite{Ru0}. While a manifold of
results were obtained thru the years for maps, 
some positive results have been established for Anosov flows  only
recently.
Notwithstanding the proof of ergodicity, and mixing, for geodesic flows
on manifolds of negative curvature \cite{Ho, AS, Si1} 
the first quantitative results consisted in the proof of exponential decay of
correlations for geodesic flows on manifolds of constant 
negative curvature in two \cite{CEG, Mo, Ra} and three \cite{Po1}
dimensions. The proof there is group theoretical in
nature and therefore ill suited to generalizations to the non
constant curvature case.\footnote{Although some partial results
for slowly varying curvature were obtained by perturbative techniques
\cite{CEG}.} The conjecture that all Axiom A
mixing flows exhibit exponential decay of correlations had already been proven false
by Ruelle \cite{Ru1, Po2} who produced piecewise constant
ceiling suspensions with arbitrarily slow rate of decay. 

The next advance was due
to Chernov \cite{Ch} who put forward the first {\em dynamical}  proof
showing sub-exponential decay of correlations for geodesic flows on
surfaces of variable negative curvature.  The basic idea was to construct a
suitable stochastic approximation of the flow (see also \cite{Li1} for
a generalization of such a point of view). 

The last substantial advance in the field
is due to the work of Dolgopyat \cite{Do1, Do2, Do3}. He
was able to use the thermodynamics formalism \cite{Si0, Ru01, Po3} and
elaborate the necessary estimate on the Perron-Frobenius operator  
to control the Laplace transform of the correlation function. As a
consequence he established exponential decay of correlations
for all Anosov flows with $\Co^{1}$ strong stable and unstable
foliations. He also gave conditions for fast decay of correlations (for
$\Co^{\infty}$ observable) in more general cases.

Unfortunately, $\Co^{1}$ strong stable and unstable
foliations seem to be a quite rare phenomenon for higher dimensional Anosov flows
\cite{PSW, Hasselblatt, SS}.
One is therefore led to think that, unless some further
geometrical structure is present, Anosov flows decay typically
slower than exponentially. 

The simplest geometrical structure that can be considered is certainly a
contact structure, geodesic flows in particular. In this case an
explicit formula by Katok and Burns \cite{KB} provides an approximation
to the {\em temporal function} which is the real quantity on which
some smoothness is required. An improvement on the error term for the above
formula, than can be found in this paper (Appendix \ref{app:Contact}, Lemma
\ref{lem:tempd}), shows that, for a Contact Anosov flow, if the strong
foliations are $\tau$-H\"older, with $\tau> \sqrt{3}-1$, then the
temporal functions is likely to be
$\Co^{1}$ (see Remark \ref{rem:smooth-fol}). On the other hand, geodesic flows
that are $a$-pinched\footnote{That is,  
such that there exists $C>0$ for which $-C\leq$ sectional curvatures
$<-aC$, clearly it must be $a\in (0,1)$. Recall that here we are
considering higher dimensional manifolds, geodesic flows on surfaces
always have $\Co^1$ foliations.} have foliations that are
$\Co^{2\sqrt a}$ (\cite{Kli} and Appendix \ref{app:Contact}; see
also \cite{HP, Ha1} for more complete results on such an issue). 
Dolgopyat result would then, at best, imply that any geodesic flow in negative
curvature which is $a$-pinched, with $a>1-\frac{\sqrt 3}2$, enjoys
exponential decay of correlations. 

Given the fact that the above numbers do not look particularly inspiring it
is then natural to guess that all Anosov Contact flows exhibit
exponential decay of correlations. This is exactly what it is proved
in the present paper (Theorem \ref{thm:main}).

To obtain such a result I built on Dolgopyat's work and on the results
in \cite{BKL} where it is introduced a functional space over which the
Perron-Frobenius operator can be studied directly, without any coding,
contrary to the previous approaches by Dolgopyat, Chernov and Pollicott.

Over such a space all the thermodynamics quantities studied by
Dolgopyat have a particularly simple analogous with a specially transparent
interpretation. It is then possible to establish a spectral gap for
the generator of the flow and this, in turn, implies exponential decay
of correlations.

The simplification of the approach is considerable as is testified by
the length of the (self-contained) proof. In addition, the transparency
of the relevant quantities allows to recognize that in certain cases
the results of Dolgopyat can be dramatically improved. To keep the
exposition as simple as possible I have chosen to restrict it to the
main case in which new results can be obtained: spectral properties of
Contact Anosov flows with respect to the Contact volume. 
This allows to choose a function space simpler than the one needed in the
general case (see \cite{BKL} for a more general choice of the Banach
space that would accommodate any Anosov flow with respect to any
equilibrium measure).

The plan of the paper is as follows. Section two starts by describing
the type of flows under consideration and the key objects used in
the proof. Then the main result is stated precisely (Theorem
\ref{thm:main}). After that a proof of the result is presented. The
proof is complete provided one assumes Lemma \ref{lem:LY}, Lemma
\ref{lem:compact} and Proposition \ref{prop:resbound}.
Lemma  \ref{lem:LY} is proven in section four. Lemma
\ref{lem:compact} is proven in section four. Section five contains the
proof of  Proposition \ref{prop:resbound} modulo and inequality, Lemma
\ref{lem:dolgo}, which is proven in section six.

Finally, for the reader convenience, the paper contains three
appendices. Appendix A contains a collection of needed--but already
well established--facts on Anosov flows. Appendix B is devoted to the
discussion of known--and less known--properties of Contact
flows. Appendix C contains few technical facts about averages that
will certainly not surprise the experts but needed to be proven somewhere.

\section{Statements and results}\label{sec:statments}

We will consider a $\Co^{4}$, $2d+1$ dimensional, connected compact
Riemannian manifold $\M$ and a
$\Co^{4}$
flow\footnote{That is $T_0=\Id$ and $T_{t+s}=T_t\circ T_s$ for each
$t,s\in\R$.} $T_t:\M\to\M$ defined on it which satisfies the following
conditions.
\begin{cond}\label{cond:tre}
At each point $x\in\M$ there exists a splitting of the tangent space
$\T_x\M=E^s(x)\oplus E^c(x)\oplus E^u(x)$. The splitting is invariant
with respect to $T_t$, $E^c$ is one dimensional and coincides with the
flow direction, in addition there exists $A,\mu>0$ such that
\[
\begin{split}
\|dT_tv\|&\leq A e^{-\mu t}\|v\|\quad \text{for each }v\in E^s\text{ and
}t\geq 0\\
\|dT_tv\|&\geq A e^{\mu t}\|v\|\quad \text{for each }v\in E^u\text{ and
}t\leq 0.
\end{split}
\]
That is, the flow is \emph{Anosov}.
\end{cond}

\begin{cond}\label{cond:due}
There exists a $\Co^{2}$ one form $\alpha$ on $\M$, such that
$\alpha\wedge(d\alpha)^d$ is nowhere zero, which  is left invariant
by $T_t$ (that is $\alpha(dT_t v)=\alpha(v)$ for each $t\in\R$ and
tangent vector $v\in\T\M$). In other words $T_t$ is a \emph{Contact Flow}.
\end{cond}

\begin{rem}
From now on I will assume $\M$ to be a Riemannian manifold with the
Riemannian volume  being the same as the contact volume
$\alpha\wedge(d\alpha)^d$. This is not
really necessary, yet it is convenient and can be done without loss of
generality.
\end{rem} 
 
With a slight abuse of notation let us define on $\Co^1(\M,\C)$ the
following group of operators 
\begin{equation}\label{eq:operator}
T_t\vf:=\vf\circ T_t\; ;\quad \Lp_t f:=f\circ T_{-t}
\end{equation}
The operator $\Lp_t$ specifies the evolution of the densities and
therefore should determine the statistical properties of the system.
Unfortunately, the spectral properties of $\Lp_t$ on $\Co^1(\M,\C)$ are not well
connected to the statistical properties of the map. To establish
such a connection it is necessary to enlarge the space. In order to do so
we must define weaker norms. Clearly such norms will need to have a
relation with the dynamical properties of the system. 

The simplest way to embed the dynamics of a system into the topology
is to introduce a \emph{dynamical distance}. In our case several
natural possibilities are available: for each $\sigma\in\R$ let
\begin{equation}\label{eq:dyndist1}
d^+_{\sigma}(x,y):=\int_0^\infty e^{\sigma t}d(T_tx,T_ty)\,dt;\quad
d^-_{\sigma}(x,y):=\int_{-\infty}^0 e^{-\sigma t}d(T_tx,T_ty)\, dt,
\end{equation}
where $d(\cdot,\cdot)$ is the Riemannian metric of $\M$.
\begin{rem}
Note that $d^+_{\sigma}$ and $d^-_{\sigma}$ are distances only if
$\sigma$ is sufficiently small (that is, negative and larger, in
absolute value, than the absolute values of all the Lyapunov
exponents), otherwise they are only
\emph{pseudo--distances}.\footnote{That is, they can attain the value
$+\infty$.} 
\end{rem}

In the present article we are interested only to the special cases of
\eqref{eq:dyndist1} considered in the following Lemma (the trivial
proof is left to the reader).
\begin{lem}\label{lem:triv}
Choose $\lambda\in(0,\mu)$ and let $d_s:=d^+_{\lambda}$ and
$d_u:=d^-_{\lambda}$. Then $d_u$ is a 
pseudo-distance on $\M$ and $d_u(T_{-t}x,T_{-t}y)\leq e^{-\lambda
t}d_u(x,y)$. In addition, $d_u$, restricted to any strong-unstable
manifold, is a smooth function and it is equivalent to the restriction
of the Riemannian metric, while points belonging to different unstable
manifolds are at an infinite distance. The analogous properties hold for $d_s$.
\end{lem} 

We can now start to describe the spaces on which we will consider the
operators $T_t$ and $\Lp_t$. First of all let us fix $\delta>0$ that will
need to be sufficiently small (how small will be specified later in
the paper) and define
\begin{equation}
\label{eq:holderunst}
H_{s,\beta}(\vf):=\sup_{d_s(x,y)\leq
\delta}\frac{|\vf(x)-\vf(y)|}{d_s(x,y)^\beta};\;\;
|\vf|_{s,\beta}:=|\vf|_\infty+H_{s,\beta}(\vf).  
\end{equation}

\begin{defin}
In the following by the Banach space $\Co^{\beta}_s(\M,\C)\subset
\Co^{0}(\M,\C)$  we will mean the closure of  $\Co^{1}(\M,\C)$ with
respect to the norm $|\cdot|_{s,\beta}$.  
Similar definitions hold with respect to the metric $d_u$ and the
Riemannian metric $d$ (giving the space of H\"older function $\Co^{\beta}$).
\end{defin}

Let us also define the unit ball
$\D_\beta:=\{\vf\in\Co_s^{\beta}(\M,\C)\;|\;|\vf|_{s,\beta} \leq 1\}$.
For a given $\beta<1$, and $f\in\Co^{1}(\M,\C)$, let  
\begin{equation}\label{eq:norms}
\begin{array}l
\|f\|_w:=\sup\limits_{\vf\in\D_1}\int_\M \vf f\\
\|f\|:=\|f\|_s+\|f\|_u\; ;\quad
\|f\|_s:=\sup\limits_{\vf\in\D_\beta}\int_\M \vf f\; ;\quad
\|f\|_u:=H_{u,\beta}(f) .
\end{array}
\end{equation}

Let $\B(\M,\C)$ and $\B_w(\M,\C)$ be the completion of $\Co^1(\M,\C)$
with respect to the norms $\|\cdot\|$ and $\|\cdot\|_w$
respectively.  Note that such spaces are separable by construction and
are all contained in $(\Co^{\beta})^*$, the dual of the
$\beta$-H\"older functions.

It is well known that the strong stable and unstable foliations for an
Anosov flow are $\tau$-H\"older (see Appendices \ref{app:Anosov},
\ref{app:Contact}  for quantitative estimates of $\tau$ and Remark
\ref{rem:tause} for the use of $\tau$ in this paper). Moreover the 
Jacobian of the holonomies associated to the stable and unstable
foliations are $\tau$-H\"older. From now on
we will assume\footnote{The square is needed only in Lemma
\ref{lem:s-vfbound}. In fact, employing the 
strategy used in \cite{BKL}, section 3.6, and refining Lemma
\ref{lem:tempd}, it may be possible to replace $\tau^2$
by $\tau$. I do not pursue this possibility since it would
complicate the proofs without any substantial addition to the present
results.} 
\begin{equation}\label{codn:beta}
\beta<\tau^2.
\end{equation}

%\penalty-2000
The {\bf main result} of the paper is the following.
\begin{thm}\label{thm:main}
For a $\Co^{4}$ Anosov Contact flow $T_t$ satisfying Conditions \ref{cond:tre} and
\ref{cond:due} the operators $\Lp_t$ form a strongly continuous group
on $\B(\M,\C)$.\footnote{In fact the only place in which the
$\Co^{4}$ hypothesis is used is in the estimate \eqref{eq:jac}. With
a bit more work, adopting the alternative approach used in \cite{BKL}
Sub-lemma 3.1.3, it is possible to reduce the needed smoothness to
$\Co^{3}$, possibly $\Co^{2+\alpha}$, but to reduce it further
some new ideas seem to be needed.}
In addition, there exists $\sigma, C_1>0$ such that,
for each $f\in\Co^{1}$, $\int f=0$, the following holds true
\[
\|\Lp_t f\|\leq C_1e^{-\sigma t}| f|_{\Co^{1}}.
\]
\end{thm} 

Clearly the above theorem implies exponential decay of correlations
for $\Co^{1}$ function:
\[
\int f\vf\circ T_t=\int\Lp_t\left[f-\int f\right]\vf+\int f\int \vf \Lp_t1=
\int f \int \vf+{\Or}(e^{-\sigma
t}|f|_{\Co^{1}}|\vf|_{s,\beta}).
\] 

In fact, a standard approximation argument
extends the result to all H\"older functions.
\begin{cor}
For each $\alpha\in(0,1)$ there exists $C_\alpha>0$ such that, for each
$f,\vf\in\Co^{\alpha}$,
\[
\left|\int f\vf\circ T_t-\int f\int \vf\right|\leq C_\alpha
|f|_{\Co^{\alpha}}|\vf|_{\Co^{\alpha}}e^{-\frac{\alpha\sigma}{2-\alpha} t}.
\]
\end{cor}

\begin{rem}
Note that Theorem \ref{thm:main} does not imply that $\Lp_1$ is a
quasicompact operator neither that it enjoys a spectral gap. This
is a reflection of the impossibility, with the ideas at hand, to
investigate directly the time one map and indicates that the result
must be pursued in a more roundabout way.
\end{rem}
The proof of Theorem \ref{thm:main} is achieved via a careful study of the
spectral properties of the generator of the group. The first step
consists in the following result proven in section \ref{sec:LY}.

\begin{lem}\label{lem:LY}
The operators $\Lp_t$ extend to a group of bounded operators on
$\B(\M,\C)$ and $\B_w(\M,\C)$; they form a strongly
continuous group. In addition, for each $\beta'<\beta$ there exists a
constant $B\geq 0$ such that, for each $f\in \B_w(\M,\C)$, $t\geq 0$,
\[
\|\Lp_t f\|_w\leq \|f\|_w
\]
and, for each $f\in \B(\M,\C)$, $t\geq 0$,
\[
\|\Lp_t f\|\leq \|f\|;\quad \|\Lp_t f\|\leq 3e^{-\lambda\beta' t}\|f\|+B\|f\|_w.
\]
\end{lem}
\relax From now on let $\beta'$ be fixed.

Accordingly the spectral radius of $\Lp_t$, $t\geq 0$, is bounded by
one. In addition, it is possible to define the generator $X$ of the 
group. Clearly, the domain $D(X)\supset\Co^{2}(\M,\C)$ and restricted to
$\Co^{2}(\M,\C)$  it is nothing else but the action of the vector
field defining the flow. 

The spectral properties of the generator depend on
the resolvent $R(z)=(z\Id-X)^{-1}$. It is well known
(e.g. see \cite{Davies}) that for all 
$z\in\C$, $\Re (z)>0$, the following holds 
\begin{equation}
\label{eq:resolvent}
R(z) f=\int_0^\infty e^{-zt}\Lp_t f dt.
\end{equation}

Thanks to \eqref{eq:resolvent} it is possible to obtain the analogue of
Lemma \ref{lem:LY} for the resolvent.

\begin{lem}\label{lem:trivialbound0}
For each $z\in\C$, $\Re (z)=a>0$, holds
\[
\|R(z)\|_w\leq a^{-1}\; ;\quad\|R(z)\|\leq a^{-1}\; ;\quad\|R(z)^nf\|\leq \frac
3{(a+\lambda\beta')^n}\|f\|+a^{-n}B\|f\|_w. 
\]
\end{lem}
\begin{proof}
The first two inequalities follow directly from formula
\eqref{eq:resolvent} and the first two inequalities of Lemma \ref{lem:LY}: 
\[
\|R(z)f\|\leq\int_0^\infty e^{-at}\|\Lp_t f\|\, dt\leq a^{-1}\|f\|.
\]
 By induction one easily obtains the formula
\begin{equation}\label{eq:respowers}
R(z)^nf=\frac 1{(n-1)!}\int_{\R^+} t^{n-1} e^{-zt}\Lp_t f dt.
\end{equation}
Using again Lemma \ref{lem:LY}
\[
\|R(z)^nf\|\leq\frac 1{(n-1)!}\int_0^\infty t^{n-1}
e^{-at}(3e^{-\lambda\beta' t}\|f\|+B\|f\|_w)
\leq \frac {3\|f\|}{(a+\lambda\beta')^{n}}+a^{-n}B\|f\|_w.
\]
\end{proof}
The next basic result (proven in section \ref{sec:quasi-comp}) is a
compactness property for the operators $R(z)$. 
\begin{lem}\label{lem:compact}
For each $a=\Re(z)>0$ the operator $R(z)$, seen as an operator
from $\B(\M,\C)$ to $\B_w(\M,\C)$, is compact.
\end{lem}

\begin{prop}
\label{lem:quasicompact}
For each $a=\Re(z)>0$ the operator $R(z)$, seen as an operator
on $\B(\M,\C)$, is quasi compact, has spectral radius $a^{-1}$ and
essential spectral radius bounded by $(a+\lambda\beta')^{-1}$.
\end{prop}
\begin{proof}
The bound on the spectral radius of $R(z)$ follows trivially from the
second inequality of Lemma \ref{lem:trivialbound0}.
While, by the third inequality of Lemma \ref{lem:trivialbound0},
Lemma \ref{lem:compact} and the usual Hennion's argument 
\cite{He} based on Nussbaum's formula \cite{Nu},  it follows that the essential
spectral radius is bounded by $(a+\lambda\beta')^{-1}$. Let us recall
the argument.
Nussbaum's formula asserts that if $r_n$ is the inf of the $r$
such that $\{R(z)^n f\}_{\|f\|\leq 1}$ can be covered by a finite
number of balls of radius $r$, then the essential spectral radius of
$R(z)$ is given by $\liminf_{n\to\infty}\sqrt[n]{r_n}$.
Let $B_1:=\{f\in\B\;|\;\|f\|\leq 1\}$. By Lemma \ref{lem:compact},
$R(z)B_1$ is relatively compact in $\B_w$. Thus, for each $\epsilon>0$ there are
$f_1,\dots,f_{N_\epsilon}\in R(z)B_1$ such that
$R(z)B_1\subseteq\bigcup_{i=1}^{N_\epsilon} 
U_\epsilon(f_i)$, where $U_\epsilon(f_i)=\{f\in\B\;|\;\|f-f_i\|_w<\epsilon\}$. 
For $f\in R(z)B_1\cap U_\epsilon(f_i)$, Lemma~\ref{lem:trivialbound0} implies that
\begin{displaymath}
\|R(z)^{n-1}(f-f_i)\|\leq\frac 3{(a+\lambda\beta')^{n-1}}\,\|f-f_i\|
                       +\frac{B}{a^{n-1}}\|f-f_i\|_w
  \leq a^{-n+1}\left\{\frac 3{(1+\lambda\beta' a^{-1})^{n-1}}+B\epsilon\right\}\ .
\end{displaymath}
Choosing $\epsilon=(1+\lambda\beta' a^{-1})^{-n+1}$ we can conclude that
for each $n\in\N$  the set $R(z)^n(B_1)$ can be covered by a finite number of 
$\|\cdot\|$--balls of radius $(3+B)(a+\lambda\beta')^{-n+1}$. 
\end{proof}

\relax For each $\zeta\in\R^+$ let $U_\zeta:=\{z\in\C\;|\;
\Re(z)>-\zeta\}$. Proposition \ref{lem:quasicompact} implies the
following corollary.\footnote{This is the equivalent of the statement
that the Laplace transform of the correlation function can be extend
to a meromorphic function in a neighborhood of the imaginary axes,
see \cite{Po3}.}
\begin{cor}
\label{cor:dispectrum}
The spectrum $\sigma(X)$ of the generator is contained in the left
half plane. The set $\sigma(X)\cap U_{\lambda\beta'}$ consists of, at most,
countably many isolated points of point
spectrum with finite multiplicity. Zero is the
only eigenvalue on the imaginary axis and has multiplicity one.
\end{cor}
\begin{proof}
If $F_z(w):=z-w^{-1}$, then $\sigma(X)=F_z(\sigma(R(z)))$. Thus the
essential spectrum of $X$ must lie outside 
$\bigcup_{\Re(z)>0}\{w\in\C\;|\; |z-w|\leq a+\lambda\beta'\}$. This is
exactly $U_{\lambda\beta'}$.  

Since $\Lp_t1=1$, and the space
$V_0:=\overline{\{f\in\Co^1(\M,\C);|\;\int f=0\}}^{\B(\M,\C)}$ 
is invariant, it follows $\sigma(X)=\{0\}\cup\sigma(X|_{V_0})$.
Next, suppose $Xf=ibf$ for some $b\in
\R$ and $f\in V_0$, $f\neq 0$, then $R(z)f=(z+ib)^{-1}f$, thus for $z=a-ib$
holds (see equation \eqref{eq:testlyu})
\[
\|f\|_u\leq
\frac{|z+ib|}{a+\beta\lambda}\|f\|_u=\frac{a}{a+\beta\lambda}\|f\|_u,
\]
that is $\|f\|_u=0$. Let $\{f_n\}\subset \Co^{1}$ be an approximating
sequence for $f$, $\vf\in\D_\beta$, and $t\in\R^+$,
\[
\left|\int f\vf\right|=\left|e^{-ibt}\int fT_t\vf\right|\leq
\left|\int f_nT_t\vf\right|+\|f-f_n\|. 
\]
Contact Anosov flows are mixing (see Corollary
\ref{cor:mix}), hence $\lim_{t\to\infty}\int f_nT_t\vf=0$.
The arbitrariness of $t$ and $n$ implies then $\int
f\vf=0$, that is $\|f\|_s=0$, which implies the contradiction $f\equiv
0$. 
\end{proof}
The above result, although rather interesting, does not suffice to
investigate the statistical properties of the system, to do so it is
necessary to exclude the presence of spectrum near the imaginary axis
(apart from $0$). This follows form the next result proven in sections
\ref{sec:resolvent}, \ref{sec:dolgo}. 
\begin{prop}
\label{prop:resbound}
There exists $b_*>0$, $\bar c >1$ and $\nu\in(0,1)$ such that for each
$z=a+ib$, $a\in[\bar c^{-1},\bar c]$, $|b|\geq b_*$, the spectral
radius of $R(z)$ is bounded by $\nu a^{-1}$.
More precisely, there exists $c^*>0$ such that, for $\bar n=\lceil c^*\ln
|b|\rceil$,
\[
\|R(z)^{\bar n}\|\leq \left(\frac \nu a\right)^{\bar n}.
\]
\end{prop}

\begin{cor}
\label{cor:spectrum}
The exists $\zeta_1<0$  such that
$\sigma(X)\cap U_{\zeta_1}=\{0\}$.
\end{cor}
\begin{proof}
By the same argument at the beginning of Corollary \ref{cor:dispectrum}, setting
$\zeta_0=\min\{\lambda\beta', \nu^{-1}-1\}$, 
$U_{\zeta_0}\cap \sigma(X)\subset\{z\in\C\;|\;
\Re(z)\in[-\zeta_0,0],\,|\Im(z)|\leq b_*\}$. By Corollary
\ref{cor:dispectrum} it follows that $U_{\zeta_0}\cap \sigma(X)$ contains
only finitely many points, from this the result follows.
\end{proof}

To conclude we need to transfer the knowledge gained on the spectrum of
$X$ into an estimate on the behavior of the 
semigroup. A typical way to do so would be to use the Weak Spectral Mapping
Theorem (\cite{Nagel}, page 91) stating that, for all $t\in\R$,
$\sigma(T_t)=\overline{\text{exp}(t\sigma(X))}$,
provided the semigroup is polynomially bounded for all times.
Unfortunately, our semigroup grows exponentially in the past.
Thus we need to argue directly. For this purpose a silly preliminary fact
is needed. 
\begin{lem}
\label{lem:smoothbound}
For each $z\in\rho(X)$ (the resolvent set) and $f\in D(X^2)$  the
following holds true
\[
\|R(z)f- z^{-1} f-z^{-2}X f\|\leq |z|^{-2}\|R(z)\|\,\|X^2 f\|.
\]
\end{lem}
\begin{proof}
This follows from the identity
$R(z)f =z^{-1} f+z^{-2}X f+z^{-2}R(z)X^2 f$, for all $f\in D(X^2)$.
\end{proof}
Next notice that, for each $a>0$ and
$f\in D(X^2)\cap \Co^0(\M,\C)$,\footnote{Just notice that, 
for $f\in D(X^2)$, $\|R(z)f\|_\infty\leq |z|^{-1}(\|X^2f\|+\|Xf\|+\|f\|)$
(see Lemma \ref{lem:smoothbound}). Hence for each $x\in\M$, $a>0$,
$R(a+ib)f(x)$ is in $L^2$ as a function of $b$. This means that for
$f\in D(X^2)$ and $x\in\M$ one can apply the inverse Laplace
transform formula and obtain the formula (\ref{eq:invLaplace})
point wise. Note that this implies only that the limit in
(\ref{eq:invLaplace}) takes place in the $L^2([0,\infty], e^{-at}dt)$
sense as a function of 
$t$. On the other hand $\Lp_t f$ is a continuous function of $t$ and,
again by Lemma \ref{lem:smoothbound}, 
$R(a+ib)f-\frac 1{a+ib}f$ is in $L^1(\R,\B)$, as a function of $b$. From this 
it follows that the limit in (\ref{eq:invLaplace}) converges in the
$\B$ norm for each $t\in\R^+$.}
\begin{equation}
\label{eq:invLaplace}
\Lp_tf= \frac 1{2\pi }\lim_{w\to \infty}\int_{-w}^{w} db\,
e^{at+ibt} R(a+ib)f.
\end{equation}

We can now conclude the section with the proof of Theorem \ref{thm:main}.

\begin{proof}[\bf Proof of Theorem \ref{thm:main}]
Let $\nu_1=\max\{\nu,\frac{4c^*}{3+4c^*}\}$ and $3\omega=\min\{\zeta_1,
(\nu_1^{-1}-1)\bar c\}$.\footnote{The constants $\nu,c^*,\bar c$ are
defined in Proposition \ref{prop:resbound}, $\zeta_1$ is defined in
Corollary \ref{cor:spectrum}.}
First of all by equation (\ref{eq:testlyu}) it follows that
\begin{equation}\label{eq:unormc}
\|\Lp_t f\|_u\leq e^{-\lambda\beta t}\|f\|_u,
\end{equation}
so we need only worry about the stable part of the norm.

Since $\int f=0$, Corollary \ref{cor:spectrum} implies that the
function $R(z)f$ is analytic in the domain $\{\Re(z)\geq -\zeta_1\}$. Then
$M:=\sup_{a\in[-2\omega,0];\,|b|\leq b_*}\|R(a+ib)f\|<\infty$,
moreover, for $a\in[-2\omega,0]$ and $|b|\geq b_*$, it follows  
\[
R(a+ib)=\left[\Id+(a-\bar c)R(\bar c+ib)\right]^{-1}R(\bar c+ib).
\] 
To see that the above formula is well defined consider that, by
hypothesis and Lemma \ref{lem:trivialbound0},  
\[
\|(a-\bar c)R(\bar c+ib)\|\leq
(1+\frac{|a|}{\bar c})\leq 1/3+2/3\nu_1^{-1}.
\]
In addition, for $\bar n=\lceil c^*\ln |b|\rceil$ Proposition
\ref{prop:resbound} implies 
\[
\|(a-\bar c)^{\bar n}R(\bar c+ib)^{\bar
n}\|\leq\left[\nu_1\left(1+\frac{|a|}{\bar c}\right)\right]^{\bar n}\leq
\left[\frac 23 +\frac{\nu_1}{3}\right]^{\bar n}.
\]
Accordingly,
\[
\begin{split}
\|\left[\Id+(a-\bar c)R(\bar c+ib)\right]^{-1}\|&\leq\sum_{n=0}^\infty
\|(a-\bar c)^nR(\bar c+ib)^n\|\\
&\leq\sum_{k=0}^\infty\|[(a-\bar c)^{\bar
n}R(\bar c+ib)^{\bar n}]^k\|\sum_{j=0}^{\bar n-1}\|[(a-\bar c)R(\bar
c+ib)]^j\|\\
&\leq\frac 9{2(1-\nu_1)^2}|b|^{c^*\ln[\frac 13+\frac 2{3\nu_1}]}
\leq \frac 9{2(1-\nu_1)^2} |b|^{1/2}.
\end{split}
\] 
Thus there exists $M_1>0$ such that, for $a\in[-2\omega,0]$ and $b\in\R$,
\begin{equation}\label{eq:upperbound}
\|R(a+ib)\|\leq M_1\sqrt{|b|}+M.
\end{equation}

To conclude we use (\ref{eq:invLaplace}) and shift
the contour of integration. For each  $f\in D(X^2)\cap\Co^0$,
\[
\Lp_t f=\frac 1{2\pi i}\int_{-2\omega+i\R}dz\, e^{zt}
R(z)f=\frac 1{2\pi i}\int_{-2\omega+i\R}dz\, e^{zt}
\left(R(z)-\frac 1z\right)f.
\]
By using Lemma \ref{lem:smoothbound} and \eqref{eq:upperbound} we have
that for each $\vf\in\D_\beta$ and $f\in D(X^2)\cap\Co^0$ holds
\[
\begin{split}
\left|\int_\M\Lp_t f\vf\right|&\leq\frac
1{2\pi}\int_{\R}db\left\|R(-2\omega+ib)f-\frac
1{-2\omega+ib}f\right\| e^{-2\omega t}\\
&\leq C\left\{\|X^2f\|+\|Xf\|+\|f\|\right\}e^{-2\omega t}.
\end{split}
\]

We have thus completed the proof for all $f\in D(X^2)\cap\Co^0$; to obtain
the announced result for $f\in\Co^{1}$ it suffices a standard
approximation argument.
Let $\phi:\R^+\to\R^+$ be a $\Co^{\infty}$ function such that
$\text{supp}(\phi)\subset (0,1)$ and $\int \phi=1$. For each $\ve>0$
define $\phi_\ve(t):=\ve^{-1}\phi(\ve^{-1}t)$ and, for each
$f\in\B(\M,\C)$,
\[
f_\ve:=\int_0^\infty \phi_\ve(t)\Lp_t f.
\]
Clearly $f_\ve\in D(X^n)\cap \Co^1$ for each $n\in\N$. More to the point
\[
\|X^2 f_\ve\|\leq \int|\phi_\ve''(t)|\|\Lp_t f\|\leq
\ve^{-2}|\phi''|_{L^1}|f|_{\Co^{1}}.
\]
In addition, if $f\in\Co^{1}$,
\[
\|f_\ve-f\|\leq \int\phi_\ve(t)|f\circ T_{-t}-f|_{\Co^{\beta}}\leq
\ve^{1-\beta}|f|_{\Co^{1}}\sup_{t\in[0,1]}|T_{-t}|_{\Co^{1}}.
\]
Accordingly, for each $f\in\Co^{1}(\M,\C)$, $\int f=0$, we have
\[
\|\Lp_t f\|\leq \|\Lp_t f_\ve\|+\|f-f_\ve\|\leq C_1e^{-2\omega
t}\ve^{-2}|f|_{\Co^{1}}+C_2\ve^{1-\beta}|f|_{\Co^{1}},
\]
and the wanted results follows by choosing
$\ve=e^{-2\omega(3-\beta)^{-1}t}$, hence $\sigma=2\omega(1-\beta)(3-\beta)^{-1}$. 
\end{proof}

\section{Proofs: Lasota--Yorke inequality}\label{sec:LY}
\begin{proof}[\bf Proof of Lemma \ref{lem:LY}]
By Lemma \ref{lem:triv}, for each $\alpha\in(0,1]$ 
\begin{equation}\label{eq:lys}
|T_t\vf|_\infty=|\vf|_\infty;\;\;
H_{s,\alpha}(T_t\vf)\leq e^{-\lambda\alpha t}H_{s,\alpha}(\vf).
\end{equation}
The first inequalities of Lemma \ref{lem:LY} are immediate since, for
$f\in\Co^1(\M,\C)$ and $\vf\in\D_\beta$ or $\vf\in\D_1$,
\[
\int_\M \vf\Lp_t\ f=\int_\M f T_t\vf.
\]
In addition, again by Lemma \ref{lem:triv}
\begin{equation}\label{eq:testlyu}
\|\Lp_t f\|_u=H_{u,\beta}(\Lp_t f)\leq e^{-\beta\lambda t}H_{u,\beta}(f)=
e^{-\lambda \beta t}\|f\|_u.  
\end{equation}

To conclude the argument we need the averaging operator\footnote{By
$W^s_\delta(x)$ we mean a ball of radius $\delta$, centered at $x$, with
respect to the metric obtained by restricting the Riemannian metric to
$W^s(x)$. By $m^s$ we designate the corresponding volume form.}
\begin{equation}\label{eq:average}
\A^s_\delta\vf(x):=\frac
1{m^s(W^s_\delta(x))}\int_{W^s_\delta(x)}\vf(z)m^s(dz) .
\end{equation}
The basic properties of such an operator consist in the following
\begin{sublem}\label{lem:average}
There exists $C>0$ such that for each $\vf\in\D_\beta$ one has
\[
\begin{array}l
|\A^s_\delta \vf-\vf|_\infty\leq C\delta^\beta |\vf|_{s,\beta}\\
H_{s,\beta}(\A^s_\delta \vf-\vf)\leq (2+C\delta)
H_{s,\beta}(\vf)+C\delta^{1-\beta}|\vf|_\infty \\
H_{s,1}(\A^s_\delta \vf)\leq C\delta^{-1}|\vf|_\infty
\end{array}
\]
\end{sublem}
The above Sub-Lemma is hardly surprising, yet its proof is a bit
technical and it is postponed to Appendix \ref{app:averages}.
By Sub-Lemma  \ref{lem:average} it follows that, given
$\vf\in\D_\beta$ and $f\in\Co^{1}$, holds
\[
\begin{split}
\int_\M f\vf=&\int_\M f\{\vf-\A^s_\delta\vf\}+\int_\M f\A^s_\delta\vf
\leq |\vf-\A^s_\delta\vf|_{s,\beta}\|f\|_s+|\A_{\delta}\vf|_{s,1}\|f\|_w\\
\leq&(C(\delta^\beta+\delta^{1-\beta})|\vf|_\infty+(2+C\delta)
H_{s,\beta}(\vf))\|f\|_s+C\delta^{-1}\|f\|_w. 
\end{split}
\]
Accordingly, remembering \eqref{eq:lys}, for each $\vf\in\D_\beta$,
\[
\begin{split}
\int_\M \Lp_tf\vf=&\int_\M f
T_t\vf\leq(C(\delta^\beta+\delta^{1-\beta})|\vf|_\infty+(2+C\delta) 
H_{s,\beta}(\vf\circ T_t))\|f\|_s+C\delta^{-1}\|f\|_w\\
\leq &
(C(\delta^\beta+\delta^{1-\beta})|\vf|_\infty+(2+C\delta)
e^{-\lambda \beta t} 
H_{s,\beta}(\vf))\|f\|_s+C\delta^{-1}\|f\|_w.
\end{split}
\]
We start by requiring $2+C\delta\leq 3$, then let
$T_0\in\R^+$ be such that $3e^{-\lambda\beta T_0}\leq
e^{-\lambda\beta'T_0}$; at last we choose $\delta$ so that
$C(\delta^{\beta}+\delta^{1-\beta})\leq e^{-\lambda\beta'T_0}$.
Thus, for each $t\leq T_0$, 
\begin{equation}\label{eq:almostdone}
\begin{split}
\|\Lp_t f\|_s&\leq 3 e^{-\lambda\beta't}\|f\|_s+C\delta^{-1}\|f\|_w\\
\|\Lp_{T_0} f\|_s&\leq
e^{-\lambda\beta'T_0}\|f\|_s+C\delta^{-1}\|f\|_w.
\end{split}
\end{equation}
For each $t\in\R^+$ we write $t=kT_0+s$, $k\in\N$, $s\in(0,T_0)$, and
we use \eqref{eq:almostdone} iteratively to obtain
\begin{equation}\label{eq:done}
\|\Lp_t f\|_s\leq 3e^{-\lambda\beta't}\|f\|_s+B\|f\|_w
\end{equation}
with $B=C\delta^{-1}(1-e^{-\lambda\beta'T_0})^{-1}$.

The strong continuity of the group follows trivially since, for each
$f\in\Co^1(\M,\C)$,\footnote{Indeed, $|f\circ
T_{-t}-f|_\infty+H_{u,\beta}(f\circ T_{-t}-f)\to 0$ as $t\to 0$.} 
\[
\lim_{t\to 0}\|\Lp_t f-f\|=0
\]
and $\Co^1(\M,\C)$ is dense in $\B(\M,\C)$ and $\B_w(\M,\C)$ by construction.
\end{proof}

\section{Proofs: Quasi-compactness of the resolvent}\label{sec:quasi-comp}

\begin{proof}[\bf Proof of Lemma \ref{lem:compact}]
The idea is to introduce approximate operators $R_\ve(z)$ (close in
norm to $R(z)$ as operators from $\B(\M,\C)$ to $\B_w(\M,\C)$) and
then consider the following sequence of maps (for some $\tau^2\geq
\beta_*>\beta>0$) 
\begin{equation}\label{eq:embedding}
\B(\M,\C) \overset{\scriptstyle
Id}{\longmapsto}\Co^{\beta}(\M,\C)^* \overset{\scriptstyle
Id}{\hookrightarrow}\Co^{\beta_*}(\M,\C)^*\overset{\scriptstyle
R_\ve(z)}{\longmapsto}\B_w(\M,\C). 
\end{equation}

The first map is clearly continuous since for each
$\vf\in\Co^{\beta}(\M,\C)$ and $f\in\B(\M,\C)$ one has
\[
\int_\M f\vf\leq \|f\|\,|\vf|_{s,\beta}\leq \|f\|\,|\vf|_{\Co^{\beta}}
\]
and thus $\|f\|_{(\Co^{\beta})^*}\leq\|f\|$. The second
is well known to be compact. Hence it suffices to prove that the last
map is continuous and the compactness of $R_\ve(z)$ as an operator from
$\B$ to $\B_w$ immediately follows. Let us postpone the proof of this
fact to Lemma \ref{lem:almostcomp}.

To define the approximate operators let us introduce the averaging
operator
\begin{equation}\label{eq:uaveg}
\A^u_\ve f(x):=Z_\ve(x)\int_{W^u_\ve(x)}f(\xi)m^u(d\xi),
\end{equation}
where $Z_\ve(x)$ is determined by the equation $\A^u_\ve 1=1$. We set
$R_\ve(z):=R(z)\A_\ve^u$.

\begin{sublem}\label{lem:u-approx}
The operators $R_\ve(z)$ satisfy\footnote{By $|||\cdot|||$ we mean the
norm of an operator viewed as an operator  from $\B(\M,\C)$
to $\B_w(\M,\C)$.}
\[
|||R(z)-R_\ve(z)|||\leq C\ve^\beta.
\]
\end{sublem}
\begin{proof}
For each $f\in\Co^{1}(\M,\C)$ and $\vf\in \Co^{0}(\M,\C)$, we have
\[
\left|\int_\M \A^u_\ve f\vf-\int_\M f\vf\right|\leq |\vf|_\infty\int_\M
dx Z_\ve(x)\int_{W^u_\ve(x)}d\xi |f(\xi)-f(x)|\leq \co\ve^\beta\|f\|_u|\vf|_\infty.
\]
Accordingly, $\|\A^u_\ve f-f\|_w\leq \co\ve^\beta\|f\|$, that is
$|||\A^u_\ve-\Id|||\leq \co\ve^\beta$. From Lemma \ref{lem:trivialbound0}
it follows $||| R_\ve(z)-R(z)|||\leq \co a^{-1}\ve^\beta$.
\end{proof}

Form Sub-Lemma \ref{lem:u-approx} and the compactness of
$R_\ve(z)$ the compactness of $R(z): \B(\M,\C)\to \B_w(\M,\C)$ is
obvious since the compact operators form a closed set.
\end{proof}

In the previous Lemma we have postponed the proof of Lemma
\ref{lem:almostcomp}. Before giving such a proof some preparatory
work is needed.

\begin{defin}\label{def:star}
Given an operator $B:\B\to\B$ we define $B^*:\B^*\to\B^*$ as
usual. Notice that if $\vf\in\D_\beta\subset \B^*$ and $B^*\vf\subset
L^\infty$ then, for each $f\in\Co^{1}$, $Bf\in L^1$, one has
\begin{equation}\label{eq:adjointg}
\int Bf \vf=\int f B^*\vf.
\end{equation}
Similar definitions hold for $\B_w$ and $\D_1$.
\end{defin}

\begin{rem}
In the following we will never need to investigate the duals $\B^*$,
$\B_w^*$; it will suffice to consider elements of $\D_\beta$ and
$\D_1$. Accordingly we will always use \eqref{eq:adjointg}.
\end{rem}

Next we isolate a result needed in the present argument but useful
also in the following.

\begin{lem}\label{lem:s-vfbound}
There exists $c>0$ such that for each $\alpha\in(0,\tau^2)$,  $\vf
\in\D_1$, $z\in\C$ with $|b|=|\Im(z)|>1$ and $a=\Re(z)>0$,
$\A^{u*}_\ve R(z)^*\vf\in \Co^{\alpha}$. More precisely
\[
|\A^{u*}_\ve R(z)^*\vf|_{\Co^{\alpha}}\leq c (|b|+\ve^{-1})|\vf|_{s,1}
\]
\end{lem}
\begin{proof}
Let $f\in\Co^{1}(\M,\C)$ and $\vf\in\D_1$, then
\[
\int_\M R_\ve(z)f \vf=\int_\M f R_\ve(z)^*\vf
\]
where $R_\ve(z)^*=\A^{u*}_\ve R(z)^*$,
\begin{equation}\label{eq:adjoint}
R(z)^*\vf(x)=\int_0^\infty e^{-zt}T_t\vf(x)\, dt
\end{equation}
by the definition of $\Lp_t$. On the other hand in Appendix
\ref{app:averages} it is shown that
\begin{equation}\label{eq:Ztilde}
\A^{u*}_\ve\vf(x)=\int_{W^u_\ve(x)}\tilde Z_\ve(x,\xi)\vf(\xi) m^u(d\xi)
\end{equation}
for some appropriate $\tau$-H\"older function $\tilde Z_\ve$ (see
Lemma \ref{lem:fub}). 
Since by \eqref{eq:adjoint} 
\[ 
\frac{d}{dt}(R(z)^*\vf)\circ T_t\bigg|_{t=0}= zR(z)^*\vf-\vf
\]
it follows  that $R(z)^*\vf$ is H\"older along the strong stable
direction and differentiable along the flow direction. Let us set
$\vf_*:=R(z)^*\vf$. 

Let $x,y$ be two points on the same strong stable manifold, and let
$\Psi$ be the stable holonomy between $W^{uc}(x)$ and
$W^{uc}(y)$. According to Lemma \ref{lem:distm}, for each
$z\in W^u_{\delta}(x)$ holds $d(z,\Psi(z))\leq \co
d(x,y)^\tau$. Moreover, $|1-J\Psi(z)|\leq \co d(x,y)^\tau$.

If $\delta^{\tau^2}\geq d(x,y)^{\tau^2}\geq \ve$ then, see Lemma \ref{lem:fub},
\[
|\A^{u*}_\ve\vf_*(x)-\A^{u*}_\ve\vf_*(y)|\leq
2\bar c |\vf|_\infty d(x,y)^\alpha\ve^{-\tau^2}.
\]
Suppose instead $d(x,y)^{\tau^2}\leq \ve$. Let $\hat\Psi : W^u(x)\to
W^u(y)$ be the weak stable holonomy ($\{\hat
\Psi(\xi)\}=W^{sc}(\xi)\cap W^u(y))$. The distance along the flow
between $\hat \Psi(\xi)$ 
and $\Psi(\xi)$ is nothing else than the temporal distance
$\Delta(y,\xi)$, (see definition at the end of appendix
\ref{app:Anosov} or Figure \ref{fig:temp}). Accordingly, Lemma
\ref{lem:tempd} yields $d(\hat \Psi(\xi),\Psi(\xi))\leq C
d(x,y)^{\tau^2}$. In addition, $W^{uc}_{\ve-c\ve
d(x,y)^{\tau^2}}(y)\subset \Psi(W^{uc}_\ve(x))\subset 
W^{uc}_{\ve+c\ve d(x,y)^{\tau^2}}(y)$.\footnote{By introducing a
coordinate system in which $W^{uc}(x)$ and $W^{s}(x)$ are linear spaces
one can represent $W^{uc}(y)$ as $\{(\xi, F(\xi))\}$ where, by the
H\"older continuity of the unstable foliation and setting $U(\xi):=D_\xi
F$, one has $\|U(\xi)\|\leq c\|F(\xi)\|^\tau$ and, by the H\"older
continuity of the unstable holonomy, $\|F(\xi)\|\leq c d(x,y)^\tau$. 
Thus, setting
$\gamma(t)=(vt,F(vt))$, with $v:=z-x$, and $z':=(z,F(z))$, one can
estimate
\[
\begin{split}
\text{dist}(y,z')=&\int_0^1\|\gamma'(t)\|dt
=\int_0^1\sqrt{\langle (v,U(vt)v),g((vt,F(vt)))(v,U(vt)v)\rangle}dt\\
=&\int_0^1\sqrt{\langle (v,0),g((vt,0))(v,0)\rangle+{\Or}(d(x,z)^2d(x,y)^{\tau^2})}\; dt=d(x,z)(1+{\Or}(d(x,y)^{\tau^2}) 
\end{split}
\]
where $g$ is the matrix defining the Riemannian metric. On the other
hand one can represent $W^s(z)$ as $\{(G(\zeta),\zeta)\}$, where
$V(\zeta):=D_\zeta G$ is bounded in norm by $c\ve^\tau$. Setting
$\Psi(z)=:(a,b)=(a, F(a))=(G(b),b)$ it follows $\|b\|\leq
cd(x,y)^\tau$, hence (provided $d(x,z)\geq d(x,y)^\tau$) 
$
\text{dist}(z',\Psi(z))\leq c\; \text{dist}((a,0),z) \leq \int_0^1\|V(bt)b\|dt
\leq c\; d(x,z)^{\tau}d(x,y)^\tau \leq c\; d(x,z) d(x,y)^{\tau^2}.
$
}
This, together with the uniform transversality between the
unstable manifold and the flow direction, implies that the symmetric
difference between $W^u_\ve(y)$ and $\hat\Psi (W^u_\ve(x))$ has a
volume bounded by a $\ve^{-1}d(x,y)^\alpha$ times the volume of
$W^u_\ve(x)$. Finally, it is easy to verify that $J\hat
\Psi=J\Psi$. Hence, remembering Lemma \ref{lem:fub},
\[
\begin{split}
|\A^{u*}_\ve\vf_*(x)-\A^{u*}_\ve\vf_*(y)|&\leq \co \left\{|\vf|_\infty
d(x,y)^\alpha(\ve^{-1}+ |b|)+
\int_{W^u_\ve(x)}\tilde
Z_\ve(x,\xi)|\vf_*(\Psi(\xi))-\vf_*(\xi)|m^u(d\xi)\right\}\\ 
&\leq\co\left\{|\vf|_\infty(\ve^{-1}+|b|)+H_{s,1}(\vf)\right\}d(x,y)^\alpha.
\end{split}
\]
To conclude note that the arguments in the proof of Sub-Lemma
\ref{lem:average} hold unchanged for $\A^{u*}_\ve$ instead of $\A^s_\ve$. Accordingly,
\[
H_{u,\alpha}(\A^{u*}_\ve\vf_*)\leq
C\ve^{-1}|\vf_*|_\infty\leq C\ve^{-1}|\vf|_\infty.
\]
While a direct computation shows
\[
\left|\frac{d}{dt}(\A^{u*}_\ve\vf_*)\circ
T_t\bigg|_{t=0}\right|_\infty\leq
C\left(|\vf_*|_\infty+\left|\frac{d}{dt}(\vf_*)\circ 
T_t\bigg|_{t=0}\right|_\infty\right)\leq C|b|\,|\vf|_\infty.
\]
Since any point in a $\delta$-neighborhood of $x$ can be reached by a
path along the stable, unstable and flow direction of length less than
const.$\delta$, the Lemma follows.
\end{proof}

We are finally able to prove the continuity of the operator
$R_\ve: \Co^{\beta_*}(\M,\C)\to B_w(\M,\C)$.

\begin{lem}\label{lem:almostcomp}
For each $\ve>0$ an $z\in\C$, $\Re(z)>0$, the operators $R_\ve(z)$ are
bounded operators from $\Co^{\beta_*}(\M,\C)^*$ to
$\B_w(\M,\C)$.
\end{lem}
\begin{proof}
By Lemma \ref{lem:s-vfbound} it follows that, for each $f\in\Co^{1}$ and
$\vf\in\D_1$, 
\[
\int_\M R_\ve(z)f \vf\leq |f|_{(\Co^{\beta_*})^*}\;
|R_\ve(z)^*\vf|_{\Co^{\beta_*}}\leq
C(|z|+\ve^{-1})|\vf|_{s,1}|f|_{(\Co^{\beta_*})^*} 
\]
which means $\| R_\ve(z)f\|_w\leq C(|z|+\ve^{-1})|f|_{(\Co^{\beta_*})^*}$
and the required result follows by an obvious density argument.
\end{proof}

\section{Proofs: Resolvent bound for large $\Im (z)$}
\label{sec:resolvent}

\begin{proof}[\bf Proof of Proposition \ref{prop:resbound}]
Lemma \ref{lem:trivialbound0} states that, for each $m,n\in\N$
and $f\in\Co^{1}(\M,\C)$, 
\begin{equation}\label{eq:starting}
\begin{split}
\|R(z)^{n+m} f\|&\leq\frac
3{(a+\lambda\beta')^m}\|R(z)^nf\|+a^{-m}B\|R(z)^nf\|_w\\
&\leq \frac 3{(a+\lambda\beta')^m a^n}\|f\|+a^{-m}B\|R(z)^nf\|_w
\end{split}
\end{equation}
hence all we need is to estimate more
precisely the weak norm of $R(z)^n f$. 

Remembering \eqref{eq:uaveg}
\begin{equation}\label{eq:begin}
\int f\vf =\int\A^u_\delta f
\vf+{\Or}(\delta^\beta \|f\|_u|\vf|_\infty)
=\int f\A^{u*}_\delta\vf+{\Or}(\|f\|_u|\vf|_\infty).
\end{equation}

Thus, for each $k,l\in\N$, $k+l=n$, and $\vf\in\D_1$ holds, by
equation \eqref{eq:begin},
\[
\int_\M R(z)^{n}f \vf=\int_\M  R(z)^kf
R(z)^{*l}\vf=\int_\M  R(z)^kf
\A^{u*}_\delta R(z)^{*l}\vf + a^{-l}{\Or}(\|R(z)^kf\|_u).
\]
To continue let 
\[
\Phi_l(\vf):=\A^{u*}_\delta R(z)^{*l}\vf .
\]
Thus, taking into account \eqref{eq:respowers} and \eqref{eq:unormc},
\begin{equation}\label{eq:due}
\int_\M R(z)^{n}f \vf=\int_\M  R(z)^kf\Phi_l(\vf)+
a^{-n}(1+a^{-1}\lambda\beta)^{-k}{\Or}(\|f\|_u).
\end{equation}

\begin{lem}\label{lem:stabcont}
There exists $c>0$ such that, for each $l\in\N$ and $\vf\in\D_1$,
\[
H_{s,\beta}(\Phi_l(\vf))\leq c |b| a^{-l}|\vf|_{s,1}.
\]
\end{lem}
\begin{proof}
The proof follows immediately from Lemma \ref{lem:s-vfbound} and
formulae \eqref{eq:respowers}, \eqref{eq:lys}.
\end{proof}
The above estimate is not particularly impressive and clearly it
can have some interest only if we can get good bounds on
$|\Phi_l(\vf)|_\infty$. This can be achieved by using an inequality
due to Dolgopyat.\footnote{Actually the original Dolgopyat estimate,
\cite{Do1}, holds for the $L^2$ norm and it is done for a different
operator in a different functional space, yet the
key cancellation mechanism due to the oscillations of the
exponential and the non joint integrability of the foliation remains
substantially identical in the two settings.} 

\begin{lem}[Dolgopyat inequality]\label{lem:dolgo}
There exists $c_*,c_1,\gamma>0$ such that, for each
$\vf\in R(z)^*(\D_1)$ and $l\geq \lceil c_*\ln |b|\rceil$, the
following holds 
\[
a^l |\Phi_l(\vf)|_\infty \leq c_1|b|^{-\gamma}l|\vf|_{s,1}.
\]
\end{lem}
The proof of the above Lemma can be found in Section \ref{sec:dolgo}.

Since equation \eqref{eq:lys} implies that, for
each $q\in\N$, $a^qR(z)^{*q}\vf\in\D_{s,1}$ and
$H_{s,\beta}(R(z)^{*q}\vf)\leq(a+\beta\lambda)^{-q} H_{s,\beta}(\vf)$ by
Lemma \ref{lem:dolgo} and Lemma \ref{lem:stabcont} it follows that
\[
|R(z)^{*k}\Phi_l(\vf)|_{s,\beta}\leq
c_4\{(1+a^{-1}\lambda\beta)^{-k}|b|+|b|^{-\gamma}l\} a^{-n}|\vf|_{s,1}.
\]
Choose $l:=\lceil c_*\ln b\rceil$, then there exists $c'>0$ and
$\nu_0\in(0,1)$ such that setting $k=\lceil c' \ln b\rceil$ equation
\eqref{eq:due} yields
\begin{equation}\label{eq:finishing}
\|R(z)^{n}f\|_w\leq c_5a^{-n}\nu_0^{n}\|f\|.
\end{equation}
The Proposition follows by \eqref{eq:starting}, \eqref{eq:finishing},
choosing $m=n=\bar n/2$ (hence $c^*=2(c_*+c')$), $\bar c=2$,
$\nu\in(\sqrt\nu_0,1)$ and $b_*$ such that $c_5(\nu_0\nu^{-2})^{n}\leq 1$. 
\end{proof}

\section{Dolgopyat Inequality}\label{sec:dolgo}

This section is devoted to the proof of Lemma \ref{lem:dolgo}. The
strategy is based on the representation \eqref{eq:respowers} (actually
on the obvious adjoint representation obtained by \eqref{eq:adjoint})
and a careful estimate of the corresponding integral.

The following simple preliminary Lemma shows that we need to worry about only
a part of the integral defining $\Phi_l(\vf)$.
\begin{lem}
There exists $\nu_*<1$ such that
\[
\left|\frac 1{(l-1)!}\int_0^{e^{-1}a^{-1}l} t^{l-1}
e^{-zt}\A^{u*}_\delta(T_t \vf)  dt\right|\leq
\nu_*^l a^{-l} |\vf|_\infty .
\]
\end{lem}
The straightforward proof is left to the reader.

Thus we can limit ourselves to consider
\[
\frac 1{(l-1)!}\int_{e^{-1}a^{-1}l}^\infty t^{l-1}
e^{-zt}\A^{u*}_\delta(T_t\vf) .
\]
To continue it is useful to localize in time. To do so we introduce a
$\Co^\infty$ function $\tp:\R\to\R$ such that $0\leq \tp\leq 1$,
$\text{supp}(\tp)\subset[-1/2,3/2]$ and with the property that
$\sum_{k=-\infty}^\infty\tp(t-k)=1$ for each $t\in \R$. Using such a
partition of unity and setting $p_0:=\lceil a^{-1}e^{-1}l\rceil$, we can write
\[
\left|\int_{0}^{\infty}t^{l-1}
e^{-zt}\A^{u*}_\delta(T_t\vf)\right|\leq\left|\sum_{k=p_0}^\infty
\int_{\R}t^{l-1}
e^{-zt}\tp(t-k)\A^{u*}_\delta(T_t\vf)\right|
+\nu_*^la^{-l}(l-1)!|\vf|_\infty .
\]

Let us analyze each of the above addenda separately.

For each $k\in\N$ holds (see \eqref{eq:Ztilde})
\[
\begin{split}
&\int_{\R}t^{l-1}
e^{-zt}\tp(t-k)\A^{u*}_\delta(T_t\vf)\\
&=\int_{\R}\tp(t-k)t^{l-1}
e^{-zt}\int_{T_k W^u_\delta(x)}\tilde Z(x,T_{-k}\xi)\vf(T_{t-k}\xi)
J_uT_{-k}(\xi),
\end{split}
\]
where by $J_uT_t$ we designate the unstable Jacobian of the map $T_t$.

To compute the above quantity it is convenient to localize
in space as well. To this end we fix a sequence of smooth partitions
of unity. There exists $c_d>0$ such that, for each $r\in(0,1)$ one can
consider a $\Co^{4}$ partition of unity
$\{\phi_{r,i}\}_{i=1}^{q(r)}$ enjoying the following
properties\footnote{It is an easy exercise to verify that partitions
with the properties below  do exist.}
\begin{enumerate}[\bf (i)]
\item for each
$i\in\{1,\dots,q(r)\}$, there exists $x_i\in\M$ such that
$\phi_{r,i}(\xi)=1$ for all $\xi\in B_r(x_i)$ (the ball of radius
$r$ centered at $x_i$) and $\phi_{r,i}(\xi)=0$ for all
$\xi\not\in B_{c_d r}(x_i)$;
\item there exists a $K>0$ such that
for each $r,i$ holds
 $\|\phi_{r,i}'(x)\|\leq
Kr^{-1}\chi_{B_{c_dr}(x_i)}(x)$;\footnote{Here, and in the following,
$\chi_A$ is the characteristic function of the set $A$.}
\item there exists $C>0$ such that $q(r)\leq Cr^{-2d-1}$. 
\end{enumerate}

Accordingly, we can write
\[
\begin{split}
&\int_{\R}t^{l-1}
e^{-zt}\tp(t-k)\A^{u*}_\delta(T_t\vf)\\
&=\sum_{i=1}{q(r)}e^{-zk}\int_{\R}\tp(t)(t+k)^{l-1}
e^{-zt}\int_{T_k W^u_\delta(x)}\phi_{r,i}(T_t\xi) \tilde
Z(x,T_{-k}\xi)\vf(T_t\xi) J_uT_{-k}(\xi) .
\end{split}
\]

From now on we will assume $b>0$, the case $b<0$ being identical.
 
\begin{figure}[ht]
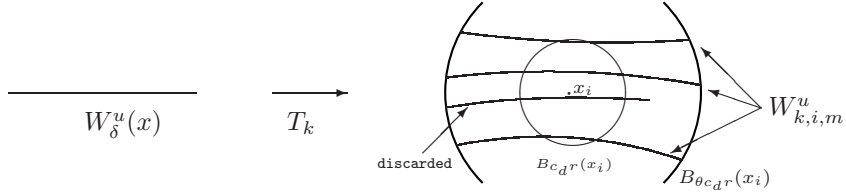
\ 
\put(0,15){\ }
\put(-50,0){\line(1,0){25}}
\put(-15,0){\vector(1,0){10}}
\put(25,0){\circle{15}}
\put(25,0){\psarc(0,0){1.7}{-45}{45}}
\put(25,0){\psarc(0,0){1.7}{135}{-135}}
\put(24,-1){$\cdot$}
\qbezier(10,8)(25,6)(40.3,7)
\qbezier(8,2)(25,4)(41.9,1)
\qbezier(8,-2)(20,0)(35,-1)
\qbezier(9.5,-7)(25,-4)(39.4,-9)
\put(-40,-5){$W^u_\delta(x)$}
\put(-13,-5){$T_k$}
\put(50,-2){\vector(-3,1){7}}
\put(50,-2){\vector(-1,1){8}}
\put(50,-2){\vector(-2,-1){12}}
\put(5,-8){\vector(1,1){6}}
\put(-1,-10){\tiny\tt discarded}
\put(51,-3){$W^u_{k,i,m}$}
\put(25,0){$\scriptstyle x_i$}
\put(20,-10){$\scriptscriptstyle B_{c_dr}(x_i)$}
\put(39,-12){$\scriptstyle B_{\theta c_dr}(x_i)$}
\put(0,-14){\ }
\caption{The manifolds $W^u_{k,i,m}$.}
\label{fig:mani}
\end{figure}

In the following we choose $\rho\in(0,\tau/8)$ and we fix
\begin{equation}
\label{eq:cond1}
r:=b^{-\varrho};\quad \varrho:=\frac{1-\tau+2\rho}{2-\tau}.
\end{equation}
It is useful to partition $T^k W^u_\delta(x)$ into submanifolds. For
each $x_i$ let us consider the connected pieces of $T^k
W^u_\delta(x)\cap B_{\theta c_dr}(x_i)$ intersecting
$B_{c_dr}(x_i)$ ($\theta$ is specified shortly). Call them
$\{W^u_{k,i,m}\}$. Among such local manifolds discard the ones such
that $\partial W^u_{k,i,m}\not\subset 
\partial B_{\theta c_dr}(x_i)$, see Figure \ref{fig:mani}.   
Clearly, if $W$ is a discarded manifold, then $T_{-k}W$ belongs to a
$\theta c_dr\lambda^{-k}$-neighborhood of $\partial W^u_\delta(x)$, hence
the total measure of the preimages of the discarded manifolds is
bounded by const.$\lambda^{-k}$. The constant 
$\theta$ is chosen so that if $\xi\in W^u_{k,i,m}\cap B_{c_dr}(x_i)$,
then $W^s_\delta(\xi)\cap W^u_{k,i,j}\neq \emptyset$, for all $j$.

Let us define $W^{uc}_{k,i,m}:=\cup_{t\in[-2,2]} T_tW^u_{k,i,m}$. 

For each $\xi\in W^{uc}_{k,i,j}$ let $t(\xi)$ be such
that $T_{t(\xi)}\xi\in  W^{u}_{k,i,j}$ and let $u(\xi):=T_{t(\xi)}\xi$.
Then
\begin{equation}\label{eq:dol1}
\begin{split}
&\int_{\R}t^{l-1}e^{-zt}\tp(t-k)\A^{u*}_\delta(T_t\vf)
=\sum_{ij}
e^{-zk}k^{l-1}\int_{W^{uc}_{k,i,j}}\tp(t(\xi))\\ 
&\;\times \left(1+\frac
{t(\xi)}k\right)^{l-1} e^{-zt(\xi)}\tilde Z(x,T_{-k}u(\xi))\vf(\xi)
\phi_{r,i}(\xi) J_uT_{-k}(u(\xi))+k^{l-1}e^{-ak}{\Or}(\lambda^{-k}|\vf|_\infty). 
\end{split}
\end{equation}

Next, for each $W^{uc}_{k,i,j}$  let
$\Psi_{k,i,j}$ be the stable holonomy between $W^{uc}_{k,i,j}$ and
$W^{uc}_{k,i,0}$. By the general theory of the holonomy maps (see Appendix
\ref{app:Anosov}) it follows that $\Psi_{k,i,j}$ is a $\tau$-H\"older
function with $\tau$-H\"older Jacobian $J\Psi_{k,i,j}$.

Notice that $T_{-k} W^u_{k,i,j}$ has size smaller than
$2c_d\lambda^{-k}r$ and thus (see Lemma \ref{lem:fub})
\[
\tilde Z(x,T_{-k}u(\xi))=Z_{k,i,j}+{\Or}(\lambda^{-k\tau}r^{\tau}).
\]
To simplify notations let us introduce the functions
\begin{equation}\label{eq:defF}
\begin{split}
F_{k,i,j}(\xi)=& \tp(t(\xi)) \left(1+\frac
{t(\xi)}k\right)^{l-1}\phi_{r,i}(\xi)
J_uT_{-k}(u(\xi)) e^{-at(\xi)}Z_{k,i,j}\\
\hat F_{k,i,j}(\xi)=&F_{k,i,j}(\Psi_{k,i,j}(\xi)) )J\Psi_{k,i,j}(\xi)
\end{split}
\end{equation}

Using the above formulae we can rewrite \eqref{eq:dol1} as
\begin{equation}\label{eq:dol2}
\begin{split}
\int_{\R}t^{l-1}\tp(t-k)
e^{-zt}\A^{u*}_\delta(T_t\vf)&=  \sum_{ij}e^{-zk}k^{l-1}\int_{W^{uc}_{k,i,0}}
e^{-ibt(\Psi_{k,i,j}(\xi))}\hat F_{k,i,j}(\xi) \vf(\xi)\\
&+k^{l-1}e^{-ak}{\Or}(\lambda^{-k\tau}|\vf|_\infty+rH_{s,1}(\vf)).
\end{split}
\end{equation}

The last preparatory step is to apply Schwartz inequality. More precisely, for each
$k,i$, we can compute
\begin{equation}\label{eq:dol3}
\begin{split}
&\left|\int_{W^{uc}_{k,i,0}} \vf(\xi)\left[\sum_j
e^{-ibt(\Psi_{k,i,j}(\xi))}\hat F_{k,i,j}(\xi)\right]\right|\leq
C|\vf|_\infty r^{(d+1)/2}\\
&\quad\quad\times\left[\sum_{j,j'}\int_{W^{uc}_{k,i,0}}
e^{-ibg_{j,j'}^0(\xi)}
\hat F_{k,i,j}(\xi) \overline{\hat F_{k,i,j'}(\xi)}\right]^{\frac 12},
\end{split}
\end{equation}
where
\begin{equation}\label{eq:defg}
g_{j,j'}^0(\xi):=t(\Psi_{k,i,j}(\xi))-t(\Psi_{k,i,j'}(\xi)) .
\end{equation}
We are finally approaching the end of the story: to conclude we must only show
that the above integral is small.

Let us perform the sum on $j'$ for each $j$. Fixed $j$ it is
convenient to express the integral on the manifold $W^{uc}_{k,i,j}$:

\[
\int_{W^{uc}_{k,i,0}}
e^{-ibg_{j,j'}^0(\xi)}
\hat F_{k,i,j}(\xi) \overline{\hat F_{k,i,j'}(\xi)}=\int_{W^{uc}_{k,i,j}}
e^{-ibg_{j,j'}(\xi)} F_{k,i,j}(\xi)\overline{F_{k,i,j,j'}(\xi)}
\]
where
\begin{equation}\label{eq:deffin}
\begin{split}
g_{j,j'}(\xi)&:=t(\xi)-t(\Psi_{k,i,j,j'}(\xi))\\
\Psi_{k,i,j,j'}(\xi)&:=\Psi_{k,i,j'}\circ\Psi_{k,i,j}^{-1}(\xi)\\
F_{k,i,j,j'}(\xi)&:=F_{k,i,j'}(\Psi_{k,i,j,j'}(\xi))
J\Psi_{k,i,j'}\circ\Psi_{k,i,j}^{-1}(\xi), 
\end{split}
\end{equation}
clearly $\Psi_{k,i,j,j'}$ is nothing else than the holonomy between
$W^{uc}_{k,i,j}$ and $W^{uc}_{k,i,j'}$.

Finally, it is convenient to divide the sum over $j'$ into two part: the sum
over nearby manifolds and the sum over manifolds at a useful
distance. Let us be more precise.

Let $y_{k,i,j}:=W^s_{2c_dr}(x_i)\cap W^{uc}_{k,i,j}$. We define the sets of indexes
$A_{k,i,j}:=\{j'\;|\; d(y_{k,i,j},y_{k,i,j'})<b^{-\varsigma}\}$ and
$B_{k,i,j}:=\{j'\;|\; d(y_{k,i,j},y_{k,i,j'})\geq b^{-\varsigma}\}$.
In the following we choose
\begin{equation}\label{eq:cond2}
\varsigma:=\frac{1-4\rho}{2-\tau}.
\end{equation}
Notice the the assumption $\rho<\tau/6$ implies that $b^{-\varsigma}$ is much
smaller than $r$, as $b$ increases.

The first step is to estimate the sum with indexes in $A_{k,i,j}$. To
do so we need the next Lemma whose proof is postponed to the end of
the section.
\begin{lem}\label{lem:mandens}
For each $\epsilon>0$, let $W^{uc}_\epsilon$ be an unstable disk of
radius $\epsilon$.
Then there exit constants $C, r_0>0$ such that for each $k\in\N$, $r_1>0$ and $x\in
\M$, calling $\{W_j\}$
the connected components of $T_kW^{uc}_\epsilon\cap B_{2r_1}(x)$, 
\[
\sum_{j\in\Omega} \sup_{\xi\in W_j}J^u_\xi T_{-k}\leq
C\,m^s(W^s_{{r_1}+\lambda^{-k}r_0}(x)), 
\] 
where $\Omega:=\{j\;|\; W_j\cap W^s_{r_1}(x)\neq \emptyset\}$.
\end{lem}

We then require
\begin{equation}\label{eq:cond3}
\lambda^{-l}r_0\leq b^{-\varsigma}.
\end{equation}
Using the above Lemma and standard distortion arguments we readily obtain
\begin{equation}\label{eq:closestuff}
\left|\sum_{j'\in A_{k,i,j}}\int_{W^{uc}_{k,i,j}}
e^{-ibg_{j,j'}}
F_{k,i,j} \overline{F_{k,i,j,j'}}\right|\leq  C J^uT_{-k}(y_{k,i,j})
b^{-d\varsigma}r^{d+1}.  
\end{equation}

We are then left with the estimate of the indexes in
$B_{k,i,j}$. To this end it is useful to make a connection with the
{\em temporal function} introduced at the end of Appendix
\ref{app:Anosov} and shown pictorially in 
Figure \ref{fig:temp}.  For each $\xi\in W^{uc}_{k,i,j}$
holds\footnote{\label{foot:fig}To apply Figure \ref{fig:temp} to the present case
set: $y=y_{k,i,j'}$, $x=y_{k,i,j}$ and $y'=T_{-t(y_{k,i,j})}u(\xi)$.}
\begin{equation}\label{eq:tempf}
g_{j,j'}(\xi)=t(\Psi_{k,i,j,j'}(\xi))-t(\xi)=
\Delta(y_{k,i,j'},T_{-t(y_{k,i,j})}u(\xi))
-t(y_{k,i,j})+t(y_{k,i,j'}).
\end{equation}

All the above work was just preparation to apply the following Lemma
(the proof can be found at the end of the section).
\begin{lem}\label{lem:clear}
For each function $G\in\Co^{\alpha}(W^{uc}_{k,i,j})$, $0<\alpha<1$,
$j'\in B_{k,i,j}$, and setting
$\bphi(u):=\phi_{i,r}(u)\phi_{i,r}(\Psi_{k,i,j,j'}(u))$, the following
holds  
\[
\left|\int_{W^{uc}_{k,i,j}}du\;e^{-ib g_{j,j'}(u)} G(u)\bphi(u)\right|\leq
C b^{-\alpha\rho}r^{d+1}|G|_{\Co^{\alpha}}.
\]
\end{lem}

Remembering \eqref{eq:defF}, \eqref{eq:dol2}, \eqref{eq:dol3},
\eqref{eq:deffin}, using
\eqref{eq:closestuff} with Lemma \ref{lem:clear}
and taking \eqref{eq:hol-bunching} into account yields\footnote{We remark
that if we choose $\alpha<\tau^2$, then
\[
\begin{split}
\sum_{j'\in
B_{k,i,j}}|\bphi^{-1}F_{k,i,j}F_{k,i,j,j'}|_{\Co^{\alpha}} 
&\leq Cl\sum_{j'\in B_{k,i,j}}|J_uT_{-k}
J_uT_{-k}\circ\Psi_{k,i,j,j'}|_{\Co^{\alpha}}  \\
&\leq Cl\sum_{j'\in B_{k,i,j}}|J_uT_{-k}J_uT_{-k}\circ\Psi_{k,i,j,j'}|_\infty
\leq Clr^dJ_uT_{-k}(y_{k,i,j}).
\end{split}
\]
}
\begin{equation}\label{eq:dol4}
\begin{split}
&\left|\int_{\R}t^{l-1}\tp(t-k)
e^{-zt}\A^{u*}_\delta(T_t\vf)(x)\right|\leq  C k^{l-1}e^{-ak}
\sum_{i}|\vf|_\infty r^{\frac{d+1}2}\\
&\quad \times\left[\sum_j
J_uT_{-k}(y_{k,i,j})\left\{r^{d+1}b^{-d\varsigma}+r^{2d+1}b^{-\alpha\rho}l
\right\}\right]^{\frac 12}+Ck^{l-1}e^{-ak} b^{-\rho}|\vf|_{s,1}\\
&\leq C k^{l-1}e^{-ak}\left\{ \sum_{i}|\vf|_\infty
r^{\frac{d+1}2}[r^{3d+1}b^{-\alpha\rho}l]^{\frac 12}+b^{-\rho}|\vf|_{s,1}\right\}\\ 
&\leq C k^{l-1}e^{-ak} b^{-\frac {\alpha\rho} 2}l^{\frac 12}|\vf|_{s,1}.
\end{split}
\end{equation}
We can finally sum over $k$ and the result follows.

We are left with the postponed proofs.

\begin{proof}[\bf Proof of Lemma \ref{lem:mandens}]
Note that the Jacobian of $T_k$ must be equal one, on the other hand
it must also be equal to the product of the stable and unstable Jacobian times a
function $\theta$ which express the ``angle'' between the stable and
unstable manifold (and hence it is H\"older). Thus, setting
$\{\xi_j\}=W_j\cap W^s_{2r_1}(x)$,
\[
\sum_{j\in\Omega} J^u_{\xi_{j}}T_{-k}=\sum_{j\in\Omega}
\theta (T_{-k}\xi_{j})  J^s_{T_{-k}\xi_{j}}T_{k}.
\]
Now, consider $T_{-k}W^s_{r_1}(x)$, clearly it will intersect $W^{uc}$ at the
points $T_{-k}\xi_j$, $j\in\Omega$. Obviously, if we consider disks
$D_j\subset T_{-k}W^s_{r_1}(x)$ centered at $T_{-k}\xi_j$ and with radius
$r_0$ sufficiently small, but depending only on $T$, they will be all
disjoint. Moreover the 
diameter of each $T_k D_j$ must be smaller than $\lambda^{-k}r_0$. This
means that $\cup_{j\in\Omega}T_k D_j$ is a collection of disjoint sets
contained in the disk  $W^s_{r_1+\lambda^{-k}r_0}(x)$. In addition, by the
usual distortion arguments, there exists $c>0$ such that
\[
J^s_{T_{-k}\xi_{j}}T_k\leq  c\, m^s(T_k D_j).
\]
Using again distortion it follows
\[
\sum_{j\in\Omega} \sup_{\xi\in W_j}J^u_\xi T_{-k}\leq
c^2|\theta|_\infty \sum_{j\in\Omega} m^s(T_k D_j)\leq C
m^s(D_{r_1+\lambda^{-k}r_0}(x)). 
\]
\end{proof}

\begin{proof}[\bf Proof of Lemma \ref{lem:clear}]
The Lemma rests on smoothness estimates for $g_{j,j'}$ which, in turn,
are obtained by estimates on $\Delta$. 
Indeed, by looking at Figure \ref{fig:temp} again it follows that for each
$\xi,\eta\in W^{uc}_{k,i,j}$ (see also footnote \ref{foot:fig})
\begin{equation}\label{eq:gdiff}
g_{j,j'}(\xi)-g_{j,j'}(\eta)=
\Delta(\Psi_{k,i,j,j'}(T_{-t(y_{k,i,j})}u(\eta)),T_{-t(y_{k,i,j})}u(\xi))
=\Delta(\Psi_{k,i,j,j'}(\eta),T_{t(\eta)-t(\xi)}\xi) .
\end{equation}
For each $y\in W^{uc}_{k,i,j}$ define $w_{j'}(y)\in E^s(y)$ by
$\text{exp}_{y}(w_{j'}(y))=\Psi_{k,i,j,j'}(y)$. Then the
normalized vectors $\hat w_{j'}(y):= 
w_{j'}(y)| w_{j'}(y)|^{-1}$  are uniformly continuous
functions. It follows that 
there exists a uniformly smooth coordinate system
$\{u_1,u_2,\dots,u_{d+1}\}:=\{u_1,\bar u)\}$ for $W^{uc}_{k,i,j}$
such that $d\alpha(\partial_{u_1},\hat w_{j'}(y))\geq c_-\|\partial_{u_1}\|$ for all
$y\in W^{uc}_{k,i,j}$. Without loss of generality we can assume
$t(u)=u_{d+1}$. Let $v(u):=\|\partial_{u_1}\|^{-1}\partial_{u_1}$.
All is needed in the following are bounds on the dependence of
$g_{jj'}$ from the coordinate $u_1$ keeping fixed the other coordinates.

For each $\bar u$, let us consider a partition $\{[a_q, a_{q+1}]\}$
of $[-2c_dr,2c_dr]$, such that 
\[
b(a_{q+1}-a_q)\,d\alpha(w_{j'}(a_q,\bar u),v((a_q,\bar u)))=2\pi.
\]
This implies
\begin{equation}\label{eq:deltaq}
2\pi c_+b^{-1}|w_{j'}(a_q,\bar u)|^{-1}\geq a_{a+1}-a_q\geq 2\pi
c_-b^{-1}|w_{j'}(a_q,\bar u)|^{-1}.
\end{equation}
Now, since $j'\in B_{k,i,j}$ it follows that, by the H\"older
continuity of the foliation,
\begin{equation}\label{eq:lastnum}
|w_{j'}(u_q)|\geq b^{-\varsigma}-C b^{-\varsigma\tau}r\geq
b^{-\varsigma}-C b^{-\frac{1-2\rho}{2-\tau}}\geq \frac 12 b^{-\varsigma}
\end{equation}
provided $b$ is large enough.
Hence, our choices imply $|a_{q+1}-a_q|<<r$ provided $b$ is
large.

Accordingly, if $u_q=(a_q,\bar u)$ and $u'=(u_1,\bar u)$, with
$u_1\in[a_q,\,a_{q+1}]$, setting $\delta_q=a_{q+1}-a_q$,  by
Lemma \ref{lem:tempd}  and \eqref{eq:gdiff} the following holds
\begin{equation}\label{eq:errorb}
|g_{j,j'}(u')-g_{j,j'}(u_q)-(u_1-a_q)d\alpha(w_{j'}(u_q),
v(u_q))|\leq C\left(|w_{j'}(u_q)|^2 \delta_q^\tau + |w_{j'}(u_q)|^\tau
\delta_q^2\right). 
\end{equation}
Indeed, $|w_{j'}(u_q)|^{\tau_-}\geq \delta_q$ and
$\delta_q^{\tau_-}\geq |w_{j'}(u_q)|$. This follows readily from
\eqref{eq:deltaq} and $\theta c_d r\geq |w_{j'}(u_q)|\geq \frac
12b^{-\varsigma}$.
By \eqref{eq:deltaq} we have $\delta_q\leq Cb^{-\frac{1-\tau+4\rho}{2-\tau}}$,
therefore \eqref{eq:cond1}, \eqref{eq:lastnum} and \eqref{eq:errorb}  yield
\begin{equation}\label{eq:der}
|g_{j,j'}(u')-g_{j,j'}(u_q)-(u_1-a_q)d\alpha(w_{j'}(u_q),
v(u_q))|\leq C b^{-1-2\rho}.
\end{equation}

Hence,\footnote{Let $m^{uc}(du)=:m(u)du$ be the measure on the manifold,
clearly $m$ is uniformly smooth.}
\[
\begin{split}
&\int_{W^{uc}_{k,i,j}}du\;m(u)e^{-ib g_{j,j'}(u)} G(u)\bphi(u)=\int d\bar
u\sum_q\int_{a_q}^{a_{q+1}} d u_1\;m(u_1,\bar u)e^{-ib g_{j,j'}(u_1,\bar u)}
G(u_1,\bar u) \bphi(u_1,\bar u)\\
& =\int d\bar
u\sum_q \bigg\{m(u_q)G(u_q)\bphi(u_q)
e^{-ibg_{j,j'}(u_q)}\int_{a_q}^{a_{q+1}} d
u_1\;e^{-ib(u_1-a_q)d\alpha(w_{jj'}(u_q), v(u_q)) }\\
&\quad\quad\quad\quad\quad+{\Or}\left((|G|_\infty
b^{-2\rho}+|G|_{\Co^{\alpha}}\delta_q^{\alpha})\int_{a_q}^{a_{q+1}}\bphi
+ |G|_{\infty}\delta_q r^{-1}\int_{a_q}^{b_q} \chi_{B_{c_d
r}(x_i)}\right)\bigg\} 
\end{split}
\]
where we have used the fact that, for each $|h|\leq \delta_q$,
$d(\Psi_{k,i,j,j'}(u_q),\,\Psi_{k,i,j,j'}(a_q+h,\bar u)) \leq
C\delta_q$ thanks to the H\"older continuity of the 
stable foliation, our choice of the parameters and since the maximal
distance between $u$ and $\Psi_{k,i,j,j'}(u)$ is bounded by a constant
times $r$.\footnote{Here is a more detailed argument: consider a
coordinate chart based at $u_q$ in which $W^{uc}_{k,i,j}$ and
$W^s(u_q)$ are linear spaces. Then $W^s(u')$, $u'=(a_q+h,\bar u)$, can
be represented as $\{(G(\xi),\xi\}_{\xi\in\R^d}$ and if
$\Psi_{k,i,j,j'}(u')=:(a,b)$, then $G(b)=a$. Setting
$d(t):=\|G(bt)\|$, by the H\"older continuity of the foliation it follows 
\[
|d'(t)|\leq C\|b\|d(t)^\tau.
\]
The above differential inequality yields $\|a\|\leq
[\delta_q^{1-\tau}+C\|b\|]^{\frac 1{1-\tau}}\leq
\delta_q[1+C\delta_q^{\tau-1}r]^{\frac 1{1-\tau}}$. The result follows since
$r<\delta_q^{1-\tau}$, the maximal ``angle'' between $W^{uc}_{k,i,j}$
and $W^{uc}_{k,i,j'}$ is bounded by $Cr^{\tau}$, and the metric in the
chart is equivalent to the Riemannian metric.}  
Continuing the above chain of inequalities yields 
\[
\begin{split}
&=\int d\bar
u\sum_q \bigg\{m(u_q)G(u_q)\bphi(u_q)e^{-ibg_{j,j'}(u_q)}\int_{0}^{\delta_q} d
u_1\;e^{-ib u_1d\alpha(w_{jj'}(u_q), v(u_q))} \bigg\}\\
&\quad+|G|_{\Co^{\alpha}} r^{d+1}{\Or}( b^{-\alpha\rho}
+b^{-\frac{2\rho}{2-\tau}})\\ 
& =\int d\bar
u\sum_q \bigg\{m(u_q)G(u_q)\bphi(u_q)e^{-ibg_{j,j'}(u_q)}\delta_q\int_{0}^{1}
ds\;e^{-2\pi i s} \bigg\}+|G|_{\Co^{\alpha}} r^{d+1}{\Or}(
b^{-\alpha\rho}) . 
\end{split}
\]
Since the inner integral equals zero exactly, the Lemma is proven. 
\end{proof}
\appendix
\section{Basic facts (Anosov flows)}\label{app:Anosov}
In this appendix we collect, for the reader's convenience, some
information on the smoothness properties of the invariant foliations
in Anosov flows that are used in the paper.

First of all, as already mentioned, for $\Co^2$ Anosov flows the
invariant distributions (sometimes called {\em splittings}) are known
to be uniformly H\"older continuous. Let us be more precise. 

For each invertible linear map $L$ let $\theta(L):=\|L^{-1}\|^{-1}$.
We define $\|d^s_xT_t\|=\|d_xT_t|_{E^s_x}\|$,
$\|d^u_xT_t\|=\|d_xT_t|_{E^u_x}\|$, $\theta(d^s_xT_t)=\theta(d_xT_t|_{E^s_x})$ and
$\theta(d^u_xT_t)=\theta(d_uT_t|_{E^u_x})$.  
Then the following holds (\cite{Hasselblatt, SS})
\begin{equation}\label{eq:dist-bunching} 
\begin{array}l
\bullet\hbox{ If there exists $\tau_d$ such that, for each
$x\in\M$, and some $t\in \R^+$,}\\
\hskip10pt\hbox{$\|d^s_xT_t\|\,\|d^u_xT_t\|^{\tau_d}<\theta(d^u_xT_t)$,
then $E^{sc}\in\Co^{\tau_d}$.}\\ 
\bullet\hbox{ If there exists $\tau_d>0$ such that, for each
$x\in\M$,  and some $t\in \R^+$,}\\
\hskip10pt\hbox{$\theta(d^s_xT_t)^{\tau_d}\,
\theta(d^u_xT_t)>\|d^s_xT_t\|^{-1}$, then $E^{uc}\in\Co^{\tau_d}$.} 
\end{array}\hfill
\end{equation}

Moreover the H\"older continuity is uniform (that is the
$\tau_d$-H\"older norm of the distributions is bounded). The above
conditions are often called {\em $\tau_d$-pinching} or {\em
bunching} conditions.

The next relevant fact is that the above splittings are
integrable. The integral manifolds are the stable and unstable
manifolds, respectively. Clearly, this implies the existence of the
weak stable and weak unstable manifolds as well. 
They form invariant continuous foliations.
Each leaf of such foliations is as smooth as the map and it is
tangent, at each point, to the corresponding distribution, \cite{KH}. In
addition, for $\Co^r$ maps, the $\Co^r$ derivatives of such manifolds
(viewed as graphs over the corresponding distributions) are uniformly
bounded, \cite{HPS}. Finally, the foliations are uniformly transversal and
$\Co^{\tau_d}$.

In the case in which both distribution are $\Co^{1+\alpha}$, it
follows by Frobenius' theorem that the Holonomy maps are
$\Co^{1+\alpha}$ (section six of \cite{PSW}). If the 
splitting is only H\"older the situation is more subtle. 

We will call {\em stable holonomy} any
holonomy constructed via the strong stable foliations and {\em unstable
holonomy} holonomies constructed by the strong unstable foliation. The basic
result on holonomies is given by the following \cite{PSW}.
\begin{equation}\label{eq:fol-bunching}
\begin{array}l
\bullet \hbox{ If there exists $\tau_h>0$ such that for some
$t\in\R^+$ and for each $x\in\M$}\\
\hskip10pt\hbox{$\|d^s_xT_t\|\,\|d^u_xT_t\|^{\tau_h}<1$, then the stable
holonomies are uniformly $\Co^{\tau_h}$.}\\
\bullet\hbox{ If there exists $\tau_h>0$ such that for some
$t\in\R^+$ and for each $x\in\M$}\\
\hskip10pt\hbox{$\theta(d^s_xT_t)^{\tau_h}\, \theta(d^u_xT_t)>1$, then the
unstable holonomies are uniformly $\Co^{\tau_h}$.}
\end{array}
\end{equation}

The relation between smoothness of holonomies and smoothness of the
foliation (in the sense that the local foliation charts are smooth) is
discussed in detail in \cite[section 6]{PSW}. Here we restrict
ourselves to what is needed in the present paper.

This is not yet enough for our purposes: we need to talk about the
smoothness of the Jacobian of the holonomies between two manifolds
$W^{uc}(x)$ and $W^{uc}(y)$.\footnote{These are a direct consequence of the
formula \cite{Ma}
\[
J\Psi(x)=\prod_{n=0}^{\infty}\frac{J^uT_{-1}(T_n\Psi(x))}{J^uT_{-1}(T_nx)}.
\]
}
\begin{equation}\label{eq:hol-bunching}
\begin{array}l
\bullet\hbox{ The stable and unstable holonomies are absolutely continuous.}\\
\bullet \hbox{ There exists $\tau_j>0$ such that }
|1-J\Psi|_\infty\leq Cd(x,y)^{\tau_j}\\
\bullet\hbox{ There exists $\tau_j>0$ such that for each $x\in\M$  the Jacobian of}\\
\hskip10pt\hbox{ the stable holonomies are uniformly $\Co^{\tau_j}$.}\\
\bullet\hbox{ There exists $\tau_j>0$ such that for each $x\in\M$ the
Jacobian of}\\
\hskip10pt\hbox{  the unstable holonomies are uniformly $\Co^{\tau_j}$.}
\end{array}
\end{equation}

The last, but not least important, object on which we need smoothness
informations is the so called {\em temporal distance}.

Fix any point $x\in\M$ and a small neighborhood
$B_\delta(x)$. Consider a smooth $2d$ dimensional manifold $\mathcal W$ containing
$W^u(x)$ and $W^s(x)$, clearly the flow is transversal to such a
manifold. On $\W$ choose a smooth coordinate system $(u,s)$
such that $\{(u,0)\}=W^u(x)$, and $\{(0,s)\}=W^s(x)$. Although it is
not necessary, for further convenience we can assume that the
coordinate system, restricted to the stable and unstable manifolds, is
the one given by the exponential map (corresponding to the metric
restricted to such manifolds). Define then
a coordinate system $(u,t,s)$ in $B_\delta(x)$ as follows:
$T_{-t}(\xi)\in\mathcal W$ and $(u,s)$ are the coordinates of
$T_{-t}(\xi)$, clearly such coordinates locally trivialize the flow.
Let $y\in B_\delta(x)\cap W^s(x)$ and $y'\in B_\delta(x)\cap
W^u(x)$. Moreover let $z'=W^u(y)\cap W^{sc}(y')$ and $z=W^s(y')\cap
W^{uc}(y)$. By construction $z$ and $z'$ are on the same flow
orbit. Thus there exists $\Delta(y,y')$ such that
$T_{\Delta(y,y')}z=z'$. The function $\Delta(y,y')$ is called {\em temporal
distance}, see Figure \ref{fig:temp} for a pictorial description.

In general the only thing that can be said is that the temporal
distance is as smooth as the strong stable and unstable foliation (see
\eqref{eq:fol-bunching}), but we will see in appendix \ref{app:Contact} that,
if some geometric structure is present, more can be said.

\section{Basic facts (Contact flows)}\label{app:Contact}
 
Given an odd dimensional (say $2d+1$) connected compact manifold $\M$ a contact
form is a $\Co^1$ 
differential $1$-form such that the $(2d+1)$-form $\alpha\wedge
(d\alpha)^d$ is non zero at every point.  

Given a flow $T_t$ on $\M$ we call it {\em contact flow} if its associated
vector field $V$ ($V(x):= \frac {d\, T_tx}{dt}\big|_{t=0}$) is such that
$d\alpha(V,v)=0$ for all vector fields $v$ and $\alpha(V)=1$, for some
contact form $\alpha$. 

Clearly the contact flow preserves the contact form and hence also the
contact volume. 

Let us start with some trivial facts showing that, for contact flows,
a bit more can be said on the quantities introduced in the previous
appendix.

\begin{lem}\label{lem:stbunst}
For a contact flow there exists a constant $C>0$ such that, for each $x\in\M$,
\[
C_0^{-1}\leq\|d_x^sT\|\theta(d_x^u T)\leq C_0\;;\quad
C_0^{-1}\leq\|d_x^uT\|\theta(d_x^s T)\leq C_0.
\]
\end{lem}
\begin{proof}
Let $v\in E^u(x)$, $|v|=1$, clearly there must exist $w\in E^s(x)$,
$|w|=1$, such that $|d\alpha(v,w)|\geq c_-|v|\,|w|$. Accordingly,
\[
c_-|v|\,|w|\leq |d\alpha(d_xT_tv,d_xT_tw)|\leq
c_+|d_xT_tv|\,|d_xT_tw|\leq c_+|d_xT_tv|\,\|d_x^sT_t\|,
\]
taking the inf on $v$ we have
\begin{equation}\label{eq:firstside}
\theta(d_x^uT_t)\,\|d_x^sT_t\|\geq c_-c_+^{-1}.
\end{equation}
On the other hand, given $w\in E^s(x)$, $|w|=1$, there must be $v\in
E^u(x)$, $|v|=1$,
such that $|d\alpha(d_xT_tw,d_xT_tv)|\geq
c_-|d_xT_tw|\,|d_xT_tv|$. Hence,
\[
c_+\geq c_-|d_xT_tw|\,|d_xT_tv|\geq c_-|d_xT_tw|\,\theta(d_x^uT_t)
\]
taking the sup over $w$ we have
\begin{equation}\label{eq:secondside}
\|d_x^sT_t\|\,\theta(d_x^uT_t)\leq c_-^{-1}c_+.
\end{equation}
The first inequality of the Lemma is then obtained putting together
\eqref{eq:firstside} and \eqref{eq:secondside}. The second inequality
follows similarly.
\end{proof}

Another trivial, but helpful, property of contact flows is the
following.
\begin{lem}\label{lem:azero}
The contact form $\alpha$ restricted to a stable or unstable manifold
must be identically zero. In addition, the form $d\alpha$ is
identically zero when restricted to a weak stable or weak unstable
manifold.
\end{lem}

\begin{proof}
The first statement is a consequence of the
invariance of $\alpha$, for example if $v$ is a stable vector then
$\alpha(v)=\lim\limits_{t\to+\infty}\alpha(dT_tv)=0$.
The second statement is proved again by invariance. Let $v,w$ be weak stable
vectors and write them as $v=v'+aV$ and $w=w'+bV$ where $v'$ and $w'$
are stable vectors. Then $d\alpha(v,w)=\lim\limits_{t\to+\infty}d\alpha(d
T_tv,d T_tw)=ab\; d\alpha(V,V)=0$.
\end{proof}
\begin{cor}\label{cor:smooth}
The distributions are smoother than indicated in Appendix
\ref{app:Anosov}: $E^u, E^s\in\Co^{\tau_d}$.
\end{cor}
\begin{proof}
Since $E^{uc}\in\Co^{\tau_d}$ and $E^u=\{v\in
E^{uc}\;|\;\alpha(v)=0\}$ the result follows trivially.
\end{proof}
\begin{rem}\label{rem:tause}
A bit more work should show that \ref{eq:fol-bunching} and
\ref{eq:hol-bunching} hold with $\tau_d$ instead than $\tau_h$ and
$\tau_j$. This it is not important for the task at hand and we will
ignore it. Throughout the paper $\tau$ will designate the best constant (less or
equal one) for which the properties in \ref{eq:dist-bunching},
\ref{eq:fol-bunching} and \ref{eq:hol-bunching} hold.
\end{rem}

The first really interesting fact concerning contact flow is given by the
following result proved in \cite{KB}, Theorem 3.6.
\begin{thm}[Katok-Burns]\label{thm:kb}
Let $\M$ be a contact manifold as above. Let $E$ be an ergodic
component of the contact flow $T$ which has positive measure and
non-zero Lyapunov exponents except in the flow direction. Then the
flow on $E$ is Bernoulli.
\end{thm}

Accordingly, by the usual Hopf argument \cite {Ho, LW}, it follows immediately.

\begin{cor}\label{cor:mix}
Let $\M$ be a connected compact contact manifold as above and let $T_t$ be an
Anosov contact flow. Then the flow is Bernoulli (and hence mixing).
\end{cor}

The proof of Theorem \ref{thm:kb} is based, among other things, on a
Lemma concerning the temporal function (see the definition at the end
of the previous appendix) which, at least for us, has an interest in
itself. Since we need it in a slightly different, stronger and more
explicit, from we will state and prove it here again.

\begin{lem}\label{lem:tempd} Assume $\alpha\in \Co^{2}$ and
conditions \eqref{eq:dist-bunching}, \eqref{eq:fol-bunching} for some $t>0$. 
Let $\bar v\in E^u(x),\,\bar  w\in E^s(x)$ be such
that $\text{exp}_x(\bar v)=y'$ and $\text{exp}_x(\bar w)=y$,\footnote{The
exponential function is with respect to the restriction of the
metric to $W^u(x)$ and $W^s(x)$, respectively.} then
\[
\Delta(y,y')=d\alpha(\bar v,\bar w)+{\Or}(\|\bar v\|^{\tau^2}\|\bar
w\|^2+\|\bar w\|^{\tau^2}\|\bar v\|^2). 
\]
In addition,
\[
\Delta(y,y')=d\alpha(\bar v,\bar w)+{\Or}(\|\bar v\|^\tau\|\bar
w\|^2+\|\bar w\|^\tau\|\bar v\|^2), 
\]
provided $\|\bar v\|^{\frac 1{\tau_-}}\leq\|\bar w\|\leq \|\bar v\|^{\tau_-}$,
$\tau_-:=\min\{\tau,(1-\tau)\}$.\footnote{The latter
limitation--although compatible with our needs--
is certainly excessive and, possibly, completely redundant. Yet, as it
will be clear from the proof, to remove it effectively it would be
necessary to have some informations on the H\"older continuity of the
foliation in $\Co^r$ topology, which seem not to be readily
available in the literature but it does hold true--at least to some extent,
see footnotes \ref{foot:dob} and \ref{foot:dob1}.}
\end{lem}
\begin{proof}
Consider the coordinate system introduced at the end of appendix
\ref{app:Anosov} to define the temporal distance. Notice that the
Euclidean metric in such coordinates gives the right measure for the
temporal distance and the distance from $x$ of points in $W^u(x)$ or
$W^s(x)$, at the same time it is uniformly equivalent to the Riemannian
metric, we can then use it without any further comment.

Let $y=(0,0,w)$ and $y'=(v,0,0)$. In coordinates the manifold $W^u(x)$
has the form $\{(u,0,0)\}$, the manifold $W^s(x)$ 
$\{(0,0,s)\}$ and the manifolds $W^{uc}(y)$, $W^{sc}(y')$ have the
form  $\{(u,t,F(u))\}$, $\{(G(s), t,s\}$, respectively. In addition,
on the one hand the smoothness of the holonomies implies $\|F\|_\infty\leq
C\|w\|^\tau$ and $\|G\|_\infty\leq C\|v\|^\tau$. On the other hand the
smoothness of the distributions implies $\|D_uF\|\leq C\|F(u)\|^\tau$
and $\|D_s G\|\leq C\|G(s)\|^\tau$. Finally, the uniform smoothness of
the manifolds implies $\|F(u)-w-D_0Fu\|\leq C\|u\|^2$,
$\|G(s)-v-D_0Gs\|\leq C\|s\|^2$.\footnote{\label{foot:dob}Actually,
here we use a very rough bound on the second derivative, one can
certainly do better. For example, since $F(u)=F(0)+D_0F(0)u+\frac
12D_0^2F(u,u)+{\Or}(\|u\|^3)\leq C\|w\|^\tau,$ it must be, at
least, $|D_0^2F|\leq C\|w\|^{\frac \tau 3}$.}

Our aim is to introduce a two dimensional manifold that captures the
essential geometric features related to $\Delta$. To do so we
introduce two smooth foliations:
$\W_u:=\{W_u(b)\;|\;b\in[0,1]\}$, $W_u(b):=\{(u,0,bF(u))\}$, and 
$\W_s:=\{W_s(a)\;|\;a\in[0,1]\}$,
$W_s(a):=\{(aG(s),0,s)\}$.\footnote{These are just a linear
interpolation between the manifolds at $x$ ant the manifolds at $y$
and $y'$, respectively.} Notice that the two above foliations are
transversal, hence for all $(a,b)\in\Sigma_0:=[0,1]^2$ is uniquely
defined the point $\{\Xi(a,b)\}:=W_u(b)\cap W_s(a)$. In fact, if we
define the function $\Phi:\R^{2d+2}\to\R^{2d}$ by
$\Phi(u,s,a,b):=(u-aG(s),s-bF(u))$, then 
$\Phi(\Xi(a,b),a,b)\equiv 0$. Since
\begin{equation}\label{eq:implic}
\frac{\partial\Phi}{\partial(u,s)}=\begin{pmatrix}\Id &-aDG\\
                                   -bDF &\Id\end{pmatrix}
=:\Id-\Lambda,
\end{equation}
where $\|\Lambda\|<1$, provided the coordinate
neighborhood has been chosen small enough. It follows that we can apply the
implicit function theorem. Accordingly $\Xi$ is a uniformly $\Co^4$
chart for the surface $\Sigma:=\Xi(\Sigma_0)$. Such a surface is
bounded by the curves $\gamma_1:= \{\Xi(a,0)\}=\{(a v,0,0)\}$,
$\gamma_2:=\{\Xi(1,b)\}$ that belongs to $W^{sc}(y')$, 
$\gamma_3:=\{\Xi(a,1)\}$ that belongs to $W^{uc}(y)$ and $\gamma_4:=
\{\Xi(0,b)\}=\{(0,0, bw)\}$. Moreover, let us set $\hat
z:=\Xi(1,1)$, clearly $\hat z$ lies on the same flow orbit of $z$ and $z'$. 
At last, consider the curves $\gamma\subset W^u(y)$ and
$\gamma'\subset W^s(y')$ obtained by transporting, along the flow
direction, $\gamma_3$ and $\gamma_2$ respectively.\footnote{The
smoothness of $W^u(y)$ and $\gamma_3$ imply trivially 
the smoothness of $\gamma$. The same considerations apply to
$\gamma'$.} Clearly $\gamma$, $\gamma_3$ and the flow line between
$\hat z$ and $z'$ bound a two dimensional manifold (contained in
$\bigcup_{t\in\R}T_t\gamma_3$), let us call it
$\Omega'\subset W^{uc}(y)$; analogously we define $\Omega$. See
Figure \ref{fig:temp} for a visual description.\footnote{Of course the
picture is a bit misleading due to a lack of dimensions. For example,
the picture does not differentiate between the 
$d$-dimensional manifold $W^u(y)$ and the curve $\gamma$.}

\begin{figure}[ht]\ 
\put(-15,0){
\put(-40,0){\line(1,0){100}}
\put(-20,-20){\line(0,1){60}}
\put(-30,-10){\line(1,1){45}}
{\linethickness{0.08mm}
\qbezier(30,-25)(35,25)(31,40)     %weak manifold of y'
\qbezier(30,-25)( 44 ,-15)(50 ,-5)
\qbezier(50,-5)(55,45)(51,60)
\qbezier(31,40)( 44 ,50)(51 ,60)   %end
\qbezier(-50,10)(20,8)(60,13)     %weak of y
\qbezier(-50,10)(-40,16)(-30,30)
\qbezier(-30,30)(40,28)(80,33)     
\qbezier(60,13)(70,19)(80,33)  %end
}
\qbezier[50](41,0)(42,10)(41,21)
\qbezier[50](-20,20)(10,18)(41,21)
\put(41,21){\line(1,1){7}} 
{\linethickness{0.3mm}    
\qbezier (41,0)(44,10)(44,24)      %strong of y'
\qbezier (-20,20)(10,20)(48,28)    %strong of y
\qbezier (44,24)(45,25)(48,28)
\put(-20,0){\line (0,1){20}}
\put(-20,0){\line (1,0){61}}
}
{\linethickness{0.005mm}          %flow line 
\qbezier (48,28)(49.5,29.2)(52.5,30.5)     
\qbezier (41,21)(38,18)(32.6,10.5)
}
\put(-19,-3){$x$}
\put(40,-4){$y'$}
\put(-23,20){$y$}
\put(39,21){$\hat z$}
\put(48,26.5){$z'$}
\put(42,24){$z$}
\put(45, -18){$W^{sc}(y')$}
\put(-42,32){$W^{uc}(y)$}
\put(55,5){\vector(-3,1){11}}
\put(56,4){$W^s(y')$}
\put(20,35){\vector(0,-1){10.5}}
\put(18,36){$W^u(y)$}
\put(10,17){${\relax_{\gamma_3}}$}
\put(38.6,7){${\relax_{\gamma_2}}$}
\put(10,5){$\Sigma$}
\put(41.5,17.5){${\relax_\Omega}$}
\put(25,21){${\relax_{\Omega'}}$}
\put(69,-11){\vector(-1,1){10}}
\put(70,-14){$W^u(x)$}
\put(-14,43){\vector(-1,-1){5}}
\put(-13,43){$W^s(x)$}
}
\caption{Definition of the {\em temporal function} $\Delta(y,y')$ and
related quantities} 
\label{fig:temp}
\end{figure}

We can now compute the required quantity.
Consider the closed curve $\Gamma$ following $\gamma_1$,
$\gamma'$ then going from $z$ to $z'$ along the flow direction and
finally coming back to $x$ via $\gamma$ and $\gamma_4$ (the bold path
in Figure \ref{fig:temp}) then  
\begin{equation}\label{eq:tempd1}
\int_\Gamma \alpha= \Delta(y,y').
\end{equation}
This is because $\alpha$ is identically zero when restricted to a
stable or unstable manifold (see Lemma \ref{lem:azero}).  On the other
hand 
\begin{equation}\label{eq:tempd2}
\int_\Gamma\alpha=
\int_{\partial\Sigma}\alpha+\int_{\partial\Omega}\alpha
+\int_{\partial\Omega'}\alpha=\int_\Sigma d\alpha
\end{equation}
where we have used Stokes theorem and the fact that $d\alpha$ is
identically zero when restricted to a weak stable or unstable
manifold (see Lemma \ref{lem:azero}).
To continue it is better to change coordinates.
\begin{equation}\label{eq:tempd3}
\begin{split}
\int_\Sigma d\alpha&=\int_{\Sigma_0}\Xi^*d\alpha=
\int_{\Sigma_0}d_{\Xi(a,b)}\alpha(D\Xi e_1,D\Xi e_2) da db\\
&=\int_{\Sigma_0}d_{x}\alpha(D\Xi
e_1,D\Xi e_2) da db+\Or(\|\Xi\|_\infty\|D\Xi
e_1\|_\infty\,\|D\Xi e_2\|_\infty),
\end{split}
\end{equation}
where we have used the fact that $\alpha$ is $\Co^2$.
By the implicit function theorem,
\begin{equation}\label{eq:derxi1}
D\Xi=-(\Id-\Lambda)^{-1}\frac{\partial\Phi}{\partial
(a,b)}=\sum_{k=0}^\infty \Lambda^k\begin{pmatrix}
				  G&0\\0&F
				  \end{pmatrix}
\end{equation}
Since all the following arguments are restricted to the hypersurface
$t\equiv 0$, from now on we will forget the $t$ coordinate.
Accordingly,
\begin{equation}\label{eq:Xiest}
\begin{split}
D\Xi
e_1&=(\Id-\Lambda)^{-1}(G,0)=(G,0)+\sum_{k=0}^\infty\Lambda^{2k}
\{\Lambda(G,0)+\Lambda^2(G,0)\}\\ 
&=(G,0)+b(\Id-\Lambda^2)^{-1}(aDG\,DF\,G,DF\,G)\\
&=(G,0)+b(a(\Id-abDG\,DF)^{-1}DG\,DF\,G,(\Id-abDF\,DG)^{-1}DF\,G)\\
&=:(v,0)+(\Delta_u v,\Delta_sv)=:\bar v+\Delta v\\
D\Xi e_2&=(0,F)+a((\Id-abDG\,DF)^{-1}DG\,F,
b(\Id-abDF\,DG)^{-1}DF\,DG\,F)\\
&=:(0,w)+(\Delta_u w,\Delta_sw)=:\bar w+\Delta w.
\end{split}
\end{equation}
Since $d_x\alpha$ is identically zero on the weak stable
and weak unstable manifold of $x$, we have
\begin{equation}\label{eq:daest}
\begin{split}
d_x\alpha(D\Xi e_1,D\Xi
e_2)&=d_x\alpha(\bar v, \bar w)+d_x\alpha(\bar v,\Delta w)
+d_x\alpha(\Delta v,\bar w)+d_x\alpha(\Delta v,\Delta w)\\
&=d_x\alpha(\bar v,\bar
w)+\Or((\|v\|+\|\Delta_uv\|)\|\Delta_sw\|+\|w\|
\|\Delta_uv\|+\|\Delta_sv\|\|\Delta_uw\|).
\end{split}
\end{equation}
The last needed estimate concerns the variation of the functions $F,G$.
\[
\begin{split}
\Delta G(a,b)&:=G(\Xi_s(a,b))-v=G(\Xi_s(a,b))-G(\Xi_s(a,0))=\int_0^b
DGD\Xi_s e_2=\int_0^b DGw+DG\Delta_sw\\
\Delta F(a,b)&:=F(\Xi_u(a,b))-w=F(\Xi_u(a,b))-F(\Xi_u(0,b))=\int_0^a
DFD\Xi_u e_1=\int_0^a DFv+DF\Delta_uv
\end{split}
\]
Remembering \eqref{eq:Xiest} we can estimate
\begin{equation}\label{eq:eraora}
\begin{split}
\|\Delta_u v\|&\leq \|\Delta G\|+ Cab \|DG\,DF (v+\Delta G)\|_\infty\leq C
\|\Delta G\|_\infty+ Cab \|DG\,DF v\|_\infty\\ 
\|\Delta_s v\|&\leq Cb\|DF\,G\|_\infty\leq Cb\|DF v\|_\infty+Cb\|DF\Delta G\|_\infty\\
\|\Delta_u w\|&\leq Ca\|DG w\|_\infty+Ca\|DG\Delta F\|_\infty\\
\|\Delta_s w\|&\leq C\|\Delta F\|_\infty+C ab\|DF\,DG w\|_\infty
\end{split}
\end{equation}
Therefore,
\[
\begin{split}
\|\Delta G\|_\infty&\leq \|DG\|_\infty\|w\|+C\|DG\|_\infty\|\Delta
F\|_\infty +C\|DG\|_\infty\|DF\,DG w\|_\infty\\
&\leq C\|DG\|_\infty\|w\|+C\|DG\|_\infty\|\Delta F\|_\infty\\
\|\Delta F\|_\infty&\leq  C\|DF\|_\infty\|v\|+C\|DF\|_\infty\|\Delta G\|_\infty
\end{split}
\]
Substituting the first in the second yields
\[
\|\Delta F\|_\infty\leq C\|DF\|_\infty\|v\|+
C\|DF\|_\infty\|DG\|_\infty\|w\|+C\|DF\|_\infty\|DG\|_\infty\|\Delta
F\|_\infty, 
\]
that is
\begin{equation}\label{eq:fine}
\begin{split}
\|\Delta F\|_\infty&\leq  C\|DF\|_\infty\|v\|+ C\|DF\|_\infty\|DG\|_\infty\|w\|\\
\|\Delta G\|_\infty&\leq  C\|DG\|_\infty\|w\|+ C\|DG\|_\infty\|DF\|_\infty\|v\|
\end{split}
\end{equation}
Using estimates \eqref{eq:fine} and \eqref{eq:eraora} in
\eqref{eq:daest} yields
\[
d_x\alpha(D\Xi e_1,D\Xi
e_2)=d_x\alpha(\bar v,\bar w)+\Or(\|DF\|_\infty\|v\|^2+\|DG\|_\infty\|w\|^2).
\]
Remembering that, by definition, $\Xi_u=aG\circ \Xi_s$ and
$\Xi_s=bF\circ \Xi_u$ we can use the above estimate in
\eqref{eq:tempd3}, \eqref{eq:tempd1} and finally obtain
\begin{equation}\label{eq:tempd4}
\Delta(y,y')=d_x\alpha(\bar v,\bar w)+ \Or(\|v\|^2\|w\|
+\|w\|^2\|v\| +\|DF\|_\infty\|v\|^2 +\|DG\|_\infty\|w\|^2)
\end{equation}
Since $\|DF\|_\infty\leq C\|F\|^\tau_\infty\leq C\|w\|^{\tau^2}$ and
$\|DG\|_\infty\leq 
C\|v\|^{\tau^2}$ the first inequality of the Lemma is proven. To
prove the second let us assume $\|w\|\geq \|v\|$, the other situation
being symmetric with respect to the exchange of the stable and unstable
direction (which corresponds to a time reversal). Remember that
$DF\in\Co^1$, hence $\|DF\|_\infty\leq C\|w\|^\tau+C\|\Xi_u\|_\infty$, thus
\[
\begin{split}
\|DF\|_\infty&\leq C\|w\|^\tau+C \|v\|+C\|DG\|_\infty\|w\|\\
\|DG\|_\infty&\leq C\|v\|^\tau+C \|w\|+C\|DF\|_\infty\|v\|
\end{split}
\]
which yields $\|DF\|_\infty\leq C\|w\|^\tau+C\|v\|$;  $\|DG\|_\infty\leq
C\|v\|^\tau+C\|w\|$.\footnote{\label{foot:dob1}Here again a better
knowledge of the size of the second derivative would improve the
result, see footnote \ref{foot:dob}.}
This proves the lemma provided  $\|v\|^\tau\geq \|w\|$. Clearly this
condition is less stringent as $\tau$ decreases, while such a
situation should be the worst case. Obviously the previous estimates 
must have been inefficient for ``large'' $\tau$. Indeed, it is possible
to do a different estimate for $\Delta F$, $\Delta G$. Suppose
$\|v\|^\tau\leq \|w\|$.
\[
\frac{d \|F\circ \Xi_u\|}{da}=\frac{\langle DFD\Xi_u e_1,
F\rangle}{\|F\|} \leq \|DF\|\|D\Xi_u e_1\|\leq C\|F\|^\tau\|G\|
\]
Integrating the above differential inequality (and the analogous one
for $G$) yields
\[
\begin{split}
\left[\|w\|^{1-\tau}-C\|G\|_\infty\right]^{\frac{1}{1-\tau}}\leq
\|F\|_\infty&\leq
\left[\|w\|^{1-\tau}+C\|G\|_\infty\right]^{\frac{1}{1-\tau}}\\
\left[\|v\|^{1-\tau}-C\|F\|_\infty\right]^{\frac{1}{1-\tau}}\leq
\|G\|_\infty&\leq \left[\|v\|^{1-\tau}+C\|F\|_\infty\right]^{\frac{1}{1-\tau}}
\end{split}
\] 
Clearly the above equations imply $\|\Delta F\|\leq \|w\|^\tau\|v\|$ and $\|\Delta
G\|\leq \|v\|^\tau\|w\|$, provided $\|v\|^{1-\tau}\geq \|w\|$. This
implies again $\|DF\|_\infty\leq C(\|w\|^\tau+\|v\|)$ and $\|DG\|_\infty\leq
C(\|v\|^\tau+\|w\|)$. In addition, $\|F\|_\infty\leq C\|w\|$,
$\|G\|_\infty\leq C\|v\|$ and $\|D\Xi e_1\|_\infty\leq
C\|G\|_\infty\leq C\|v\|$, $\|D\Xi e_2\|_\infty\leq C\|w\|$. Using
such estimates in \eqref{eq:eraora}, \eqref{eq:daest} and
\eqref{eq:tempd3} yields 
\[
\Delta(y,y')=d_x\alpha(\bar v,\bar
w)+\Or(\|v\|^2\|w\|^\tau+\|w\|^2\|v\|^\tau).
\]
\end{proof}

\begin{rem}\label{rem:smooth-fol}
It may be possible to optimize Lemma \ref{lem:tempd}  by pushing
forward (or backward) the picture
until $d(T_kx,\,T_ky)=d(T_kx,\,T_ky')$; of course one would need to
be rather careful by properly estimating distortion. At any rate, the best result
one can hope for is that if $\tau>\sqrt{3}-1$, then
$\Delta(y,y')=d\alpha(v,w)+o(|v|)$. That is, $\Delta$ is differentiable
with respect to $y'$ and the derivative is $\Co^{\tau}$. We do not
push matters in such a direction since it is not necessary for the
purpose at hand.
\end{rem}

\section{Averages}\label{app:averages}

We start with a long overdue proof.
\begin{proof}[\bf Proof of Sub-Lemma \ref{lem:average}]

Clearly
\begin{equation}\label{eq:easy}
|\A^s_\delta\vf|_\infty\leq |\vf|_\infty;\quad
|\A^s_\delta\vf-\vf|_\infty\leq \delta^{\beta}H_{s,\beta}(\vf).
\end{equation}
The estimate of the smoothness of $\A^s_\delta\vf$ is a bit more
subtle, to investigate it is convenient to introduce an appropriate
coordinate system.

Since all the quantities are related to the same stable manifold, form
now on we will consider the Riemannian metric restricted to the
stable manifold.

Given $x,y$ belonging to the same stable manifold, we first identify the
tangent spaces at $x$ and $y$ by parallel transport, then we consider
normal coordinates at $x$ and at $y$. Clearly in such
coordinates the balls $W^s_\delta(x)$ and $W^s_\delta(y)$ are actual
balls of radius $\delta$; of course this it is not the case for
$W^s_\delta(y)$ in the normal coordinates at $x$. We call
$I_{xy}:\T_x\M\to\T_y\M$ the isometry that identifies the tangent spaces
and we define the map $\Upsilon_{xy}:\M\to\M$ as
\[
\Upsilon_{xy}(z)=\text{exp}_y(I_{xy}\text{exp}^{-1}_x(z)).
\]
where $\text{exp}$ is the exponential map defined by the metric on the
stable manifold. 

First of all notice that, by construction
\begin{equation}\label{eq:ball}
\Upsilon_{xy}(W^s_\delta(x))=W^s_\delta(y).
\end{equation} 
Next, to study $\Upsilon_{xy}$ we describe it in the normal
coordinates of the point $x$, we will then identify all the tangent
spaces by the Cartesian structure of such a chart. Calling, as usual,
$\Gamma^k_{ij}$ the Christoffel symbols, the 
equation of parallel transport for a vector $v$ along the curve
$\gamma$ reads
\[
\frac{d v^k}{dt}=-\sum_{ij}\Gamma^k_{ij}v^i\frac{d\gamma^j}{dt}.
\]
Moreover, in the normal coordinates of the point $x$,\footnote{Clearly
the smoothness of the metric will depend on the smoothness of the
tangent planes (that in our case are uniformly $\Co^{3}$), see
\cite{KH}. Accordingly $\Gamma$ will be uniformly $\Co^{2}$.}
\[
|\Gamma^k_{ij}(\xi)|\leq C|\xi|.
\]

Assuming $d(x,y)\leq \delta$, we are interested only in a region
contained in the ball $W^s_{2\delta}(x)$, thus $|\Gamma^k_{ij}|_\infty\leq
C\delta$. Hence, by a standard use of Gronwald inequality,
\begin{equation}\label{eq:ixy}
|I_{x,y}v-v|\leq C_1d(x,y)^2|v|.
\end{equation}
Arguing in the same manner on the equations defining the geodesics,
and taking into account \eqref{eq:ixy}, it follows that 
\begin{equation}\label{eq:uxy}
d(\Upsilon_{xy}(z),z)\leq (1+C_2\delta) d(x,y).
\end{equation}

This implies that the symmetric difference $W^s_\delta(x)\Delta
W^s_\delta(y)$ is contained in the spherical shell $W^s_{\delta+C_2
d(x,y)}(x)\backslash W^s_{\delta-C_2 d(x,y)}(x)$ whose measure is
proportional to $\delta^{d-1}d(x,y)$.

To see that the Jacobian is close to one a bit more work is
needed. Namely we must linearize the geodesic equations along the
geodesic. This is a standard procedure and it is best done via the
Jacobi fields \cite{DoCa}. By using Gronwald again, and the fact that
the manifolds are uniformly $\Co^{4}$, yields 
\begin{equation}\label{eq:jac}
|J \Upsilon_{xy}-1|_\infty\leq C_3 d(x,y).
\end{equation}
From this it follows immediately 
\begin{equation}\label{eq:easys}
H_{s,1}(\A^s_\delta\vf)\leq C\delta^{-1}|\vf|_\infty.
\end{equation}
We can then conclude by using \eqref{eq:ball}, \eqref{eq:uxy} and
\eqref{eq:jac}, 
\begin{equation}\label{eq:balls}
\begin{split}
\left|\int_{W^s_\delta(x)}\vf-\int_{W^s_\delta(y)}\vf\right|
&\leq\int_{B_\delta(0)}|\vf(\xi)\rho(x,\xi)-
\vf\circ\Upsilon_{xy}(\xi)\rho(y,\Upsilon_{xy}(\xi))
J\Upsilon_{xy}(\xi)|\, d\xi\\
&\leq \left[(1+c_2\delta)^\beta
H_{s,\beta}(\vf)d(x,y)^\beta+C_4|\vf|_\infty
d(x,y)\right]m^s(W_\delta^s(x)).
\end{split}
\end{equation}
\end{proof}

Next we need an estimate of how much two nearby manifolds can drift
apart.

\begin{lem}\label{lem:distm}
There exists a constant $C>0$ such that for each $x\in\M$ and $y\in
W^s_\delta(x)$ holds\footnote{Here by ``dist'' we mean the Hausdorff
distance.} 
\[
\text{\rm dist}(W^u_\delta(x),W^u_\delta(y))\leq C d(x,y)^{\tau}.
\] 
\end{lem}
\begin{proof}
Clearly $d(W^u_\delta(x),\,W^u_\delta(y))$ is bounded by the distance
computed along the stable manifold. For each $\xi\in W^u_\delta(x)$
consider the unstable holonomy between $W^{sc}(x)$ and $W^{sc}(\xi)$. Let
$\{\eta\}:= W^{sc}(\xi)\cap W^u_\delta(y)$. By \ref{eq:hol-bunching} it
follows $d_s(\xi,\eta)\leq C d_s(x,y)^\tau$. From this the Lemma follows.
\end{proof}
The other needed results concerning averages are all based on a sort
of change of order of integration formula. Although such a result may
already exist in some form in the literature (after all it is a sort
of Fubini with respect to a foliation with H\"older smoothness), I
find it more convenient to derive it in the following.

To proceed it is helpful to choose special coordinates in which the
unstable, or the 
stable manifolds, are straight. Let us do the construction for the 
unstable manifold, the one for the stable being similar.

First notice that such  a straightening can be only local, we can then
choose an appropriate covering $\{U_i\}$ of $\M$ (appropriate means
that the open sets must be sufficiently small) and perform the wanted
construction in each open set $U_i$. 

Let $U$ be a sufficiently small open ball. Let us choose a coordinate
system in $U$, since the Euclidean norm in the coordinate is
equivalent to the Riemannian length we will use it instead and we
will, from now on, confuse $U$ with its coordinate
representation.

It is particularly convenient to choose the chart in such a way that, given a
preferred point $\bar x\in U$, $\{(u,0)\}_{u\in\R^{d_u}}=W^u(\bar x)$ and
$\{(0,s)\}_{s\in\R^{d_s+1}}=W^{sc}(\bar x)$.

At this point we can define the function
$H:\R^{d_u+d_s+1}\to\R^{d_s+1}$ by the requirement
\[
\{(u,\,H(u,s))\}_{u\in\R^{d_u}}=W^u((0,s)).
\]
Clearly this implies $H(0,s)=s;\quad H(u,0)=0$.
We define then the change of coordinates
\[
\Psi(\bar u,\bar s)=(\bar u,\, H(\bar u,\bar s))
\]
in the coordinates $(\bar u,\bar s)$ the unstable manifolds are just
all the vector spaces of the type $\{(\bar u, a)\}$ for some
$a\in\R^{d_s+1}$.

In addition, a trivial computation shows that, calling $J\Psi$ the
Jacobian of the change of coordinates $\Psi$, we have that $J\Psi(\bar
u,\bar  s)$ is nothing else than the Jacobian of the unstable holonomy
between $\{(0,\xi)\}_{\xi\in\R^{d_s+1}}$ and $\{(\bar u,\xi)\}_{\xi\in\R^{d_s+1}}$.

\begin{lem}\label{lem:fub}
There exists $\bar c>0$ such that the kernel $\tilde Z_\ve$, defined in
\eqref{eq:Ztilde}, satisfies
\[
|\tilde Z_\ve|_{\Co^{\tau}}\leq \bar c |Z_\ve|_\infty,
\]
moreover $\tilde Z(x,\xi)$ is Lipschitz with respect to the second
variable, limited to the flow direction, with Lipschitz constant $\bar
c |Z_\ve|_\infty$. 
\end{lem}
\begin{proof}
Since all the relevant quantities are local quantities, we can compute
in a chart $\Psi$ as above.

Let $U$ be an open set in the chart and consider $f:\M^2\to\C$
supported in $U^2$. Then
\[
\int_\M m(dx) \int_{W^u_\delta(x)}f(x,\xi)m^u(d\xi)=
\int_{\{(x,\xi)\in U^2\;|\;d^u(x,\xi)\leq \delta\}}f(x,\xi)m^u(d\xi)m(dx).
\]
Now we set $\Xi_\delta:=\{(x,\xi)\in U^2\;|\;d^u(x,\xi)\leq \delta\}$ and
we change variables: $x=\Psi(u,s)$ and $\xi=\Psi(u',s)$.
\[
\int_\M m(dx)
\int_{W^u_\delta(x)}f(x,\xi)m^u(d\xi)=\int_{\Xi_\delta}
f(\Psi(u,s),\Psi(u',s)) \rho(u_1,s)J\Psi(u,s) du'\,du\,ds,
\]
where $\rho\circ\Psi^{-1}$ is a uniformly $\tau$-H\"older function.
Accordingly, 
\[
\tilde
Z_\ve(x,\xi)=\frac{J\Psi(x)\rho(\Psi^{-1}(\xi))}{J\Psi(\xi)\rho(\Psi^{-1}(x))}
Z_\ve(x).
\]
The smoothness of $\tilde Z_\ve(x,\xi)$ follows then from previous results on
holonomy smoothness and the smoothness of $Z_\ve$. In turn, the latter
is proven exactly as in equation \eqref{eq:balls} exchanging the 
r\^ole of the stable and unstable manifolds and setting $\vf=1$.
\end{proof}


\begin{thebibliography}{99}
\footnotesize
\bibitem{AS}  D.V.Anosov, Ya.G. Sinai,  {\em Certain smooth ergodic
systems}, Russian Math. Surveys 22 (1967), no. 5, 103--167.
\bibitem{BKL} M.Blank, G. Keller, C. Liverani, {\em
Ruelle-Perron-Frobenius spectrum for Anosov maps}, Nonlinearity, 15,
n.6, 1905-1973 (2002). 
\bibitem{Ch} N.Chernov, {\em Markov approximations and decay of correlations for
Anosov flows}, Annals of Mathematics, {\bf 147} (1998), 269--324.
\bibitem{CEG} P.Collet, H. Epstein, G. Gallavotti, {\em Perturbations of geodesic
flows on surfaces of constant negative curvature}, Communications in Mathematical
Physics, {\bf 95} (1984), 61--112.
\bibitem{Davies} E.B.Davies, {\em One-Parameter semigroups}, Academic
Press, London, (1980).
\bibitem{DoCa} M.P. Do Carmo, {\em Riemannian Geometry,} Birkh\"auser, 
Boston (1992).
\bibitem{Do1} D.Dolgopyat, {\em Decay of correlations in Anosov
flows}, Annals of mathematics, {\bf 147} (1998), 357-390.
\bibitem{Do2} D.Dolgopyat, {\em Prevalence of rapid mixing in hyperbolic 
flows}, Ergodic Theory and Dynamical Systems, {\bf 18}, 1097-1114 (1998).
\bibitem{Do3} D.Dolgopyat, {\em Prevalence of rapid mixing-II:
topological prevalence}, Ergodic Theory and Dynamical Systems, {\bf 20}
(2000), no. 4, 1045--1059.
\bibitem{Hasselblatt} B.Hasselblatt, {\em Regularity of the Anosov
splitting and of horospheric foliations}, Ergod. Th.\& Dynam.
Sys. \textbf{14} (1994), 645-666.
\bibitem{Ha1} B.Hasselblatt, {\em Horospheric foliation and relative
pinching}, J.Differential Geometry, {\bf 39}, (1994) 57-63.
\bibitem{He} H. Hennion, {\em Sur un th\'eor\`eme spectral et son
application aux noyaux Lipchitziens}, Proceedings of the American
Mathematical Society, {\bf 118} (1993), 627--634.
\bibitem{HP} M.Hirsch, C.Pugh, {\em Smoothness of horocycle
foliations},  J.Differential Geometry, {\bf 10}, (1975) 225-238.
\bibitem{HPS} M.Hirsch, C.Pugh, M.Shub, {\em Invariant Manifolds},
Lecture Notes in Math. \textbf{583} (1977).
\bibitem{Ho}  E. Hopf,  {\em Statistik der geodätischen Linien in
Mannigfaltigkeiten negativer Krümmung},
Ber. Verh. Sächs. Akad. Wiss. Leipzig 91, (1939). 261--304.  
\bibitem{KB} A.Katok, K.Burns, {\em Infinitesimal Lyapunov functions,
invariant cone families and stochastic properties of smooth dynamical
systems}, Ergodic Theory and Dynamical Systems, {\bf 14}, 757-785,
(1994).
\bibitem{KH} A.Katok, B.Hasselblatt, {\em Introduction to the Modern
Theory of Dynamical Systems}, Encyclopedia of Mathematics and its
Applications, {\bf 54}, G.-C.Rota editor, Cambridge University Press,
Cambridge (1995).
\bibitem{Kli} W.P.A.Klingenberg, {\em Riemannian Geometry}, second
edition, Walter de Gruyter, Berlin, New York (1995).
\bibitem{Li0} C.Liverani, {\em Decay of Correlations}, Annals of
Mathematics, {\bf 142} (1995), 239-301.
\bibitem{Li1} C.Liverani, {\em Flows, Random Perturbations and Rate of
Mixing}, Ergodic Theory and Dynamical Systems, {\bf 18} (1998),
1421--1446. 
\bibitem{LW}  C.Liverani,  M.Wojtkowski, {\em Ergodicity in
Hamiltonian Systems}, Dynamics Reported,  {\bf 4}, C.K.R.T. Jones,
U.Kirchgraber, H.O.Walther edts., Springer-Verlag, Berling,
Heidelgerg, New York (1995) 130-202.
\bibitem{Ma} R.Ma\~ne, {\em Ergodic Theory and Differentiable Dynamics},
Springer-Verlag, Berlin Heidelberg (1987).
\bibitem{Mo} C.Moore, {\em Exponential decay of correlation coefficients for
geodesic flows}, in {\sl Group Representation Ergodic Theory}, Operator Algebra and
Mathematical Physics, Springer, Berlin, (1987).
\bibitem{Nagel} G.Greiner, R. Nagel, {\em Spectral theory}, in {\sl
One-Parameter semigroups of Positive Operators}, Nagel edt., Lecture
Notes in Mathematics, {\bf 1184}, Springer-Verlag, Berlin, (1980).
\bibitem{Nu}  R.D. Nussbaum, {\em The radius of the essential spectrum},
Duke Math. J., {\bf 37} (1970), 473-478.
\bibitem{Po1} M.Pollicott, {\em Exponential mixing for geodesic flow on
hyperbolic three manifold}, Journal of Statistical Physics, {\bf 67} (1992),
667-673.
\bibitem{Po2} M.Pollicott, {\em A complex Ruelle-Perron-Frobenius
theorem and two counterexamples}, Eergodic Theory and Dynamical
Systems, {\bf 4} (1984), 135-146.
\bibitem{Po3} M.Pollicott, {\em On the rate of mixing of Axiom A
flows}, Inventiones Mathematicae, {\bf 85} (1985), 413--426.
\bibitem{PSW} C.Pugh, M.Shub, A.Wilkinson, {\em H\"older foliations}, 
Duke Math. J. 86 (1997), no. 3, 517--546. {\em Correction to: ``H\"older
foliations" [Duke Math. J. 86 (1997), no. 3, 517--546]},
Duke Math. J. 105 (2000), no. 1, 105--106. 
\bibitem{Ra} M. Ratner, {\em The rate of mixing for geodesic and horocycle flows},
Ergodic Theory and Dynamical Systems, {\bf 7} (1987), 267--288.
\bibitem{RS} M.Reed, B.Simon, {\em Methods of modern mathematical
physics. Fouries-Analisys, self-Adjointness}, vol. {\bf 2}, Academic
Press, New York, San Francisco, London (1972).
\bibitem{Ru0} D.Ruelle, {\em A measure associated with Axiom A
attractors}, Amer. J. Math., {\bf 98}, 616--654 (1976).
\bibitem{Ru01} D.Ruelle, {\em Thermodynamics formalism},
Addison-Wesley, New York, 1978.
\bibitem{Ru1} D.Ruelle, {\em Flots qui ne m\'elange pas
exponentialment}, C.R.Acad. Sc. Paris, {\bf 296} (1983), 191--193.
\bibitem{Si1} Ya.G.Sinai, {\em Geodesic flows on compact surfaces of
negative curvature}, Soviet Math. Dokl. 2 (1961) 106--109. 
\bibitem{Si0} Ya.G.Sinai, {\em Gibbs measures in ergodic theory},
Russian Math. Surveys {\bf 27}, 21-69 (1972).
\bibitem{SS} J.Schmeling, R.Siegmund--Schultze, {\em H\"older Continuity of the
Holonomy Map for Hyperbolic Basic Sets I}, in {\sl Ergodic Theory and Related Topics
III}, Proceedings International Conference (G\"ustrow, Germany, 1990),
Springer Lecture notes in Mathematics, eds. U.Krengel, K.Richter and
V.Warstat, Springer, Berlin, {\bf  1514} (1992), 174--191.
\end{thebibliography}
\end{document}